\documentclass{imsart}

\usepackage[autonum]{tchdr}
\usepackage[numbers,sort&compress]{natbib}

\startlocaldefs
\newcommand{\muinv}{\mu^{\leftarrow}}
\newcommand{\nuinv}{\nu^{\leftarrow}}
\endlocaldefs

\begin{document}
	
	\begin{frontmatter}
		
		\title{Truncated Simulation and Inference in\\ Edge-Exchangeable Networks}
		\thankstext{t1}{This work is supported by a National Sciences and Engineering Research Council of Canada (NSERC) 
			Discovery Grant and Discovery Launch Supplement.}
		\runtitle{Truncated Edge-Exchangeable Networks}
		
		\author{\fnms{Xinglong} \snm{Li}\ead[label=e1]{xinglong.li@stat.ubc.ca}}
		\and
		\author{\fnms{Trevor} \snm{Campbell}\ead[label=e2]{trevor@stat.ubc.ca }}
		
		\address{Department of Statistics\\ The University of British Columbia\\ email: \url{xinglong.li@stat.ubc.ca}; \url{trevor@stat.ubc.ca}}
		\runauthor{X. Li and T. Campbell}
		
		\begin{abstract}
Edge-exchangeable probabilistic network models generate edges as an
i.i.d.~sequence from a discrete measure, providing a simple means for
statistical inference of latent network properties. The measure is often
constructed using the self-product of a realization from a Bayesian
nonparametric (BNP) discrete prior; but unlike in standard BNP models, the
self-product measure prior is not conjugate the likelihood, hindering the development
of exact simulation and inference algorithms. Approximation via finite truncation of the
discrete measure is a straightforward alternative, but incurs an unknown
approximation error. 
In this paper, we develop methods for 
forward simulation and posterior inference in random self-product-measure
models based on truncation, and provide theoretical guarantees on the quality
of the results as a function of the truncation level.
The techniques we present are general and extend to the broader 
class of discrete Bayesian nonparametric models.
%
\end{abstract}

		
	\end{frontmatter}

\section{Introduction}
Probabilistic generative models have for many years been key tools in the
analysis of network data \cite{Erdos59,Holland83}. Recent work in the
area \cite{cai2016edge, williamson2016nonparametric, crane2016edge, campbell2018exchangeable, janson2018edge,
caron2017sparse, todeschini2016exchangeable, herlau2016completely, veitch2015class, borgs2017sparse, caron2017sparsity} has begun to incorporate 
the use of nonparametric discrete measures, in an effort to address the limitations of traditional models in
capturing the sparsity of real large-scale networks \cite{orbanz2014bayesian}. 
These models construct a discrete random measure $\Theta$ (often a completely random measure, or CRM \cite{kingman1967completely}) on a space $\Psi$,
associate each atom of the measure with a vertex in the network,
and then use the \emph{self-product} of the measure---i.e., the 
measure $\Theta\times\Theta$ on $\Psi^2$---to represent the magnitude of
interaction between vertices. 

While the inclusion of a nonparametric measure enables capturing sparsity,
it also makes both generative simulation and posterior inference via Markov chain Monte Carlo (MCMC) 
[\citealp{Robert04}; \citealp[Ch.~11, 12]{Gelman13}] 
 significantly more challenging.
In standard Bayesian models with discrete nonparametric measures---such as the Dirichlet process mixture model \cite{Escobar95}
or beta process latent feature model \cite{Griffiths05}---this issue is typically addressed
by exploiting the conjugacy of the (normalized) completely random measure prior 
and the likelihood to marginalize the latent infinite discrete measure \cite{broderick2018posteriors}. 
But in the case of nonparametric
network models, however, there is no such reprieve; the self-product of a completely random
measure is not a completely random measure, and exact marginalization is typically not possible.

Another option is to \emph{truncate} the discrete CRM to have finitely many atoms, and perform
simulation and inference based on the truncated CRM. Exact truncation schemes 
based on auxiliary variables \cite{Teh07,Kalli11,Zhu20} 
are limited to models where certain tail probabilities can be computed exactly.
Fixed truncation \cite{blei2006variational, blei2007correlated, wang2011online, doshi2009variational},
on the other hand, is easy to apply to any CRM-based model; but it involves an approximation with potentially unknown error. 
Although the approximation error of truncated CRM models has been thoroughly studied in past 
work \cite{campbell2019truncated,Ishwaran:2001,Ishwaran:2002,doshi2009variational,Roychowdhury:2015},
these results apply only to generative simulation---i.e., not inference---and do not apply 
to self-product CRMs that commonly appear in network models. 

In this work, we provide tractable methods for both generative simulation and posterior inference
of discrete Bayesian nonparametric models based on truncation, as well as rigorous
theoretical analyses of the error incurred by truncation in both scenarios.
In particular, our theory and methods require only the ability to 
compute bounds on---instead of exact evaluation of---intractable tail probabilities.
Our work focuses on the case of self-product-measure-based
\emph{edge-exchangeable} network sequences \cite{cai2016edge,crane2016edge,cai2015edge,crane2015edge},
whose edges are simulated \iid conditional on the discrete random product measure $\Theta\times\Theta$, 
but the ideas here apply without much effort to the broader class of discrete Bayesian nonparametric models.
As an intermediate step of possible independent interest, we also show that the nonzero rates
generated from the \emph{rejection representation} \cite{rosinski2001series} of a Poisson process
have the same distribution as the well-known but typically intractable \emph{inverse L\'evy} or \emph{Ferguson-Klass} 
representation \cite{ferguson1972representation}. This provides a novel method for simulating
the inverse L\'evy representation, which has a wide variety of uses in applications of Poisson processes \cite{wolpert1998poisson, teh2009indian, teh2010hierarchical}.

\section{Background}
\subsection{Completely random measures and self-products}
A \emph{completely random measure} (CRM) $\Theta$ on $\Psi$ is a random measure such that for any collection of $K\in\nats$ disjoint measurable sets $A_1,...,A_K \subset \Psi$, the values $\Theta(A_1), ... , \Theta(A_K)$ are independent random variables \cite{kingman1967completely}.  In this work, we
focus on discrete CRMs taking the form $\Theta =
\sum_{k}\theta_{k}\delta_{\psi_k}$, where $\delta_{x}$ is a Dirac measure on
$\Psi$ at location $x\in \Psi$ (i.e., $\delta_x(A) = 1$ if $x\in A$ and 0
otherwise), and $(\theta_k, \psi_k)_{k=1}^\infty$ are a sequence of
\emph{rates} $\theta_k$ and \emph{labels} $\psi_k$ generated from a Poisson
process on $\reals_+ \times \Psi$ with mean measure $\nu(\dee\theta) \times
L(\dee\psi)$.  Here $L$ is a diffuse probability measure, and $\nu$
is a $\sigma$-finite measure satisfying $\nu(\reals_+) = \infty$,
which guarantees that the Poisson process has countably infinitely 
many points almost surely.  The
space $\Psi$ and distribution $L$ will not affect our analysis; thus as a
shorthand, we write $\distCRM(\nu)$ for the distribution of $\Theta$: 
\[
\Theta \defined \sum_{k}\theta_k\delta_{\psi_k} \sim \distCRM(\nu).  \label{eq:crmnu}
\]

One can construct a multidimensional measure $\Theta^{(d)}$ on $\Psi^d$,
$d\in\nats$ from $\Theta$ defined in \cref{eq:crmnu} by taking its 
\emph{self-product}.  In particular, we define
\[
\Theta^{(d)} \defined \sum_{i\in\nats_{\neq}^d}\vartheta_{i}\delta_{\zeta_i}, 
& &
\vartheta_{i} \defined \prod_{j=1}^d \theta_{i_j},
& &
\zeta_i \defined (\psi_{i_1}, \psi_{i_2}, ... , \psi_{i_d}),
\label{eq:dcrm} 
\]
where $i$ is a $d$-dimensional multi-index, and $\nats_{\neq}^d$ is the set of such indices
with all distinct components. Note that
$\Theta^{(d)}$ is no longer a CRM, as it does not satisfy the independence
condition. 

\subsection{Series representations} 

To simulate a realization $\Theta\sim\distCRM(\nu)$---e.g., as a first step in 
the simulation of a self-product measure $\Theta^{(d)}$---the rates $\theta_k$ and 
labels $\psi_k$  may be generated in sequence using a \emph{series representation}~\cite{rosinski1990series} of the CRM.
In particular, we begin by simulating the 
jumps of a unit-rate
homogeneous Poisson process $(\Gamma_k)_{k=1}^\infty$ on $\reals_+$ in increasing order. 
For a given distribution $g$ on $\reals_+$ and 
nonnegative measurable function $\tau : \reals_+\times\reals_+\to\reals_+$, we set
\[ 
{\tiny }\Theta =
\sum_{k=1}^{\infty}\theta_k\delta_{\psi_k}, \quad \theta_k = 
\tau(U_k, \Gamma_k), \quad U_k\distiid g, \quad \psi_k\distiid L. \label{eq:seriesrep} 
\] 
Depending on the particular choice of $\tau$ and $g$, one can construct several different series representations of a CRM~\cite{campbell2019truncated}. 
For example, the \emph{inverse L\'evy representation} \cite{ferguson1972representation} has the form
\[
\theta_k = \nuinv(\Gamma_k), \quad \nuinv(x) \defined \inf \left\{y : \nu\left(\left[y, \infty\right)\right)\leq x\right\}.
\label{eq:levyrep}
\]
In many cases, computing $\nuinv(x)$ is intractable, making it hard to
generate $\theta_k$ in this manner. 
Alternatively, we can generate a series of rates from $\distCRM(\nu)$ with the
\emph{rejection representation}~\cite{rosinski2001series}, which has the form
\[
\theta_k =
T_k\mathds{1}\Big(\frac{\dee\nu}{\dee\mu}(T_k) \geq U_k\Big), \quad
T_k = \mu^{\leftarrow}(\Gamma_k), \quad
U_k\distiid\mathrm{Unif}[0,1],
\label{eq:rejecrep} 
\] 
where $\mu$ is a measure on $\reals_+$ chosen such that 
$\frac{\dee\nu}{\dee\mu}\leq 1$ uniformly and $\muinv(x)$ is easy to calculate
in closed-form. While there are many other sequential representations of CRMs \cite{campbell2019truncated},
the representations in \cref{eq:seriesrep,eq:levyrep,eq:rejecrep} are broadly applicable
and play a key role in our theoretical analysis.

\subsection{Edge-exchangeable graphs}\label{sec:edgeexch}

Self-product measures $\Theta^{(d)}$ of the form \cref{eq:dcrm} with $d=2$ have 
recently been used as priors in a number of probabilistic 
network models \cite{cai2016edge,williamson2016nonparametric}.\footnote{There are also network
models based on self-product measure priors that do not generate edge-exchangeable sequences \cite{caron2017sparse,veitch2015class}; it is likely that many of the techniques in the present work would extend
to these models, but we leave the study of this to future work.}
The focus of the present work are those
models that associate each $\psi_k$ with a vertex, each tuple $\zeta_i = (\psi_{i_1},\dots, \psi_{i_d})$
with a (hyper)edge, and then build a growing sequence of networks by sequentially generating edges 
from $\Theta^{(d)}$ in rounds $n=1, \dots, N$. There are a number of choices for how to construct
such a sequence.
For example, in each round $n$ we may add multiple edges via an \emph{independent likelihood process} $X_n \sim \mathrm{LP}(h, \Theta^{(d)})$ defined by
\[
X_n \defined \sum_{i\in\nats_{\neq}^d} x_{ni}\delta_{\zeta_i}, & & x_{ni}\distind h(\cdot | \vartheta_i),
\label{eq:multiedge}
\]
where $x_{ni} = k$ denotes that there were $k$ copies of edge $\zeta_i$ added at round $n$,
and $h(\cdot | \vartheta)$ is a probability distribution on $\nats \cup \{0\}$. We denote the mean
$\mu(\vartheta) \defined \sum_{k=0}^\infty k\cdot h(k\given \vartheta)$ and 
 probability of 0 under $h$ to be
$\pi(\vartheta) \defined h(0 | \vartheta)$ for convenience.
By the Slivnyak-Mecke theorem \cite{last2017lectures}, if $h$ satisfies
\[
\int_{\reals_+^d} \mu\left(\prod_{j=1}^d\theta_j\right) \prod_{j=1}^d \nu(\dee\theta_j) < \infty,  \label{eq:hcond}
\]
then finitely many edges are added to the graph in each round.
We make this assumption throughout this work.
Alternatively, if
\[
\int_{\reals_+} \min(1, \theta)\nu(\dee\theta) < \infty, \label{eq:normcond}
\]
then $\Omega\defined \Theta^{(d)}(\Psi^d) < \infty$,
and we may add only a single edge per round $n$ via a \emph{categorical likelihood process} $X_n \sim \mathrm{Categorical}(\Theta^{(d)})$
defined by
\[
X_n \defined \delta_{\zeta_{I_n}}, & & I_n \sim \mathrm{Categorical}\left( \left(\frac{\vartheta_i}{\Omega}\right)_{i\in\nats_{\neq}^d} \right). \label{eq:categorical}
\]
This construction has appeared in \cite{williamson2016nonparametric}, where $\Theta$ follows a Dirichlet process, which can be seen as a normalized gamma process \cite{Ferguson73}. 
Note that our definition of $\mathrm{Categorical}(\Theta^{(d)})$ involves simulating from the normalization; we use this definition
to avoid introducing new notation for normalized processes.

Using either likelihood process, the edges in the network after $N$ rounds are
\[
\sum_{n=1}^N X_n \defined \sum_{i\in\nats_{\neq}^d} x_i\delta_{\zeta_i}, \qquad x_i \defined \sum_{n=1}^N x_{ni},
\]
i.e., $x_i \in \nats\cup\{0\}$ represents the count of edge $i$ after $N$ rounds. 

There are three points to note about this formulation. 
First, since the atom locations $\zeta_i$ are not used, we can 
 represent the network using only its array of edge counts 
\[
E_N \defined (x_i)_{i\in\nats_{\neq}^d} \in \mcN_d,\label{eq:EN}
\]

where $\mcN_d$ denotes the set of integer arrays indexed by $\nats_{\neq}^d$ with
finitely many nonzero entries. Note that $\mcN_d$ is a countable set.
Second, by construction, the distribution
of $E_N$ is invariant to reorderings of the arrival of edges, and thus 
the network is \emph{edge-exchangeable} \cite{cai2016edge,crane2016edge,campbell2018exchangeable,janson2018edge}.
Finally, note that the network $E_N$ as formulated in \cref{eq:EN} is in general a directed multigraph
with no self-loops (due to the restriction to indices $i\in\nats_{\neq}^d$ rather than $\nats^d$).
Although the main theoretical results in this work are developed in this setting, 
we provide an additional result in \cref{sec:trunclemma}  to translate to other 
common network structures (e.g. binary undirected networks).

\section{Truncated generative simulation}\label{sec:truncation}
\begin{figure}[t!]
\begin{subfigure}{0.45\textwidth}
    \centering\includegraphics[width=0.8\textwidth]{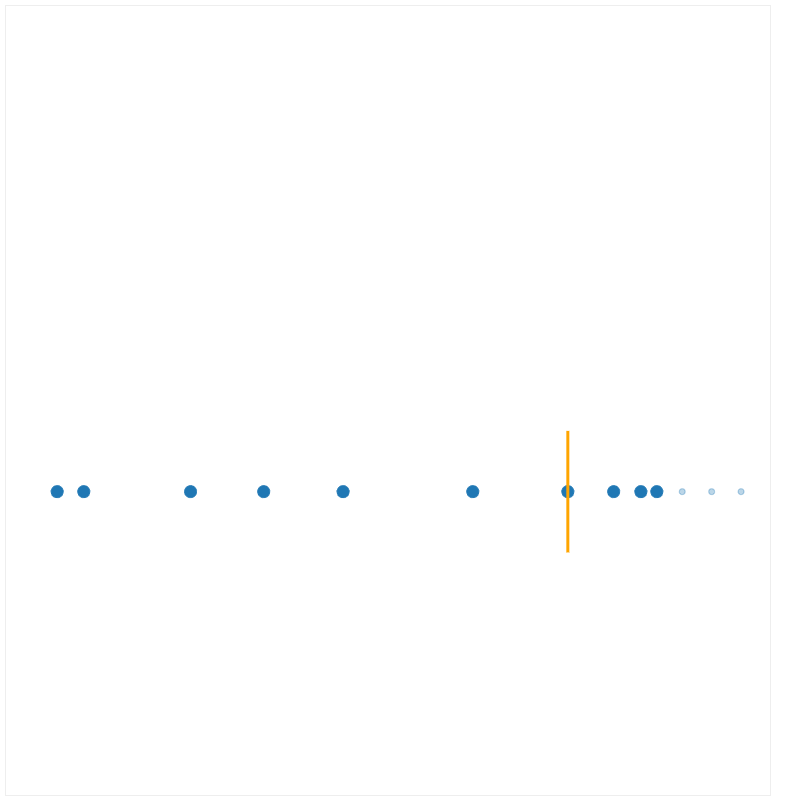}
\caption{Truncation of $\Theta$}\label{fig:trunc1}
\end{subfigure}
\begin{subfigure}{0.45\textwidth}
    \centering\includegraphics[width=0.8\textwidth]{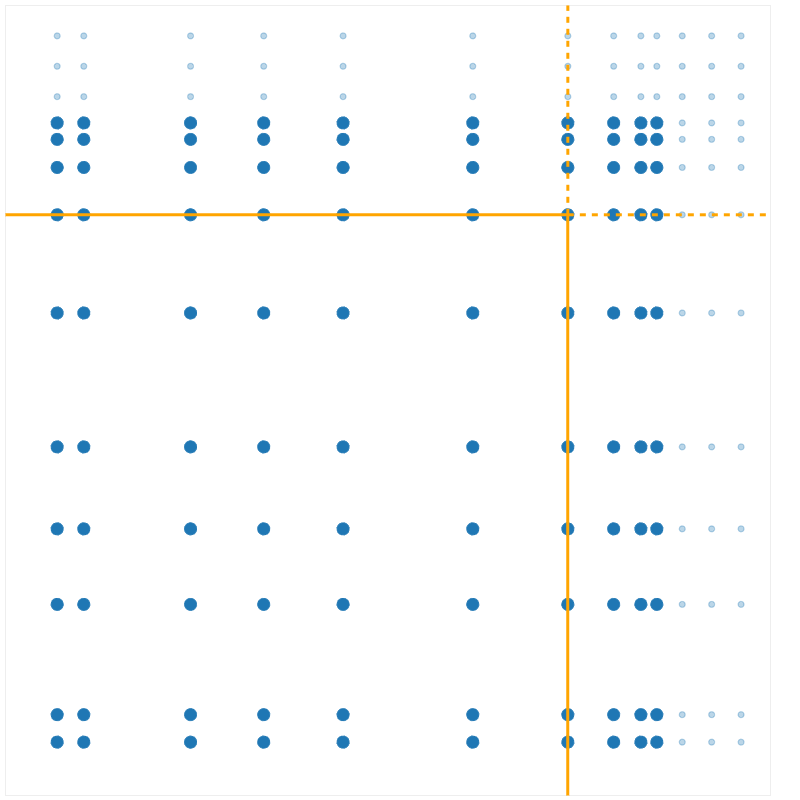}
\caption{Truncation of $\Theta^{(2)}$}\label{fig:trunc2}
\end{subfigure}
\caption{An illustration of the difference between truncation of CRMs ($d=1$) and self-product CRMs with $d=2$. Intuitively, increasing $d$ means that a higher proportion of mass is discarded in the truncation process.}\label{fig:truncdemo}
\end{figure}

In this section, we consider the tractable generative simulation of network
models via truncation, 
and analyze the approximation error incurred in doing so
as a function of $K\in\nats$ (the truncation level) and number of rounds of generation $N\in\nats$.
In particular, to construct a truncated self-product measure, we first split the underlying
CRM $\Theta$ into a truncation and tail component, 
\[
\Theta = \Theta_K + \Theta_{K+},& & \Theta_K = \sum_{k=1}^{K}\theta_k\delta_{\psi_k}, & &\Theta_{K+} = \sum_{k = K+1}^{\infty}\theta_k\delta_{\psi_k},
\]
and 
construct the self-product $\Theta_K^{(d)}$ from the truncation $\Theta_K$ as in \cref{eq:dcrm}.
\cref{fig:truncdemo} provides an illustration of the truncation of $\Theta$ and $\Theta^{(2)}$,
showing that  $\Theta^{(2)}$ can be decomposed into a sum of four parts, 
\[ \Theta^{(2)}
= \Theta^2 = (\Theta_K + \Theta_{K+})^2 = \Theta_K^{(2)} +
\left(\Theta_K\times\Theta_{K+} + \Theta_{K+}\times\Theta_{K} + \Theta_{K+}^2\right).  
\]
Thus, while we
only discard $\Theta_{K+}$ in truncating $\Theta$ to $\Theta_K$, we discard
three parts in truncating $\Theta_K^{(2)}$ to $\Theta_K^{(2)}$; and in general,
we discard $2^d - 1$ parts of $\Theta^{(d)}$ when truncating it to
$\Theta_K^{(d)}$. We therefore intuitively might expect higher truncation error when
approximating $\Theta_K^{(d)}\approx \Theta^{(d)}$ than when approximating
$\Theta_K\approx \Theta$; in \cref{sec:crmtrunc,sec:ncrmtrunc}, we will show that this is 
indeed the case.

Once the measure is truncated, a truncated network---based on $\Theta_{K}^{(d)}$---can be constructed in the same manner as 
the original network using the independent likelihood process \cref{eq:multiedge} or 
categorical likelihood process \cref{eq:categorical}.
We denote $E_{N,K} = (x_{i,K})_{i\in\nats^d_{\neq}} \in\mcN_d$ to be the corresponding edge set of the truncated network up
round $N$, where $x_{i,K} = 0$ for any index $i\in\nats^d_{\neq}$ such that some component $i_j > K$.
We keep $E_N$ and $E_{N,K}$ in the same space in order to compare their distributions in \cref{sec:trunclemma,sec:crmtrunc,sec:ncrmtrunc}.

\subsection{Truncation error bound}\label{sec:trunclemma}

We formulate the approximation error incurred by truncation as the 
total variation distance between the marginal network distributions, i.e., 
of $E_N$ and $E_{N,K}$. The first step in the analysis of truncation 
error---provided by \cref{lem:marginal}---is to show that 
this is bounded above by the probability that
there are edges in the full network $E_N$ involving vertices beyond the truncation level $K$.
To this end, we denote the \emph{maximum vertex index} of $E_N$ to be 
\[
I_N \defined \max_{i\in \nats_{\neq}^d} \left(\max_{j \in [d]} \,\, i_j\right) \quad \text{s.t.} \quad x_i > 0,
\]
and note that by definition, $I_N \le K$ if and only if all edges in $E_N$ fall in the truncated region.
\bnlem\label{lem:marginal}
Let $\Theta = \sum_{k=1}^\infty \theta_k\delta_{\psi_k}$ be a random discrete measure,
and $\Theta_K=\sum_{k=1}^K\theta_k\delta_{\psi_k}$ be its truncation to $K$ atoms. 
Let $\Theta^{(d)}$ be the self-product of $\Theta$, and $\Theta^{(d)}_K$ be the self-product
of $\Theta_K$.
Let $P_{N}$ and $P_{N,K}$ be the marginal distributions
of edge sets $E_N, E_{N,K}\in\mcN_d$ under either the independent or categorical likelihood process.
Then 
\[
\tvd{P_{N}}{P_{N,K}} \le 1- \Pr\left( I_N \leq K\right).
\]
\enlem
As mentioned in \cref{sec:edgeexch}, the network $E_N$ is in general a directed multigraph with no self loops.
However, \cref{lem:marginal}---and the downstream
truncation error bounds presented in \cref{sec:crmtrunc,sec:ncrmtrunc}---also
apply to any 
graph $E'_N=(x'_i)_{i\in\nats^d_{\neq}}$  that is a function of the original graph $E'_N = f(E_N)$
such that truncation commutes with the function, i.e., $E'_{N,K} = f(E_{N,K})$.
For example,  to obtain a truncation error bound for the common setting of undirected binary graphs, we
generate the directed multigraph $E_N$ as above and construct the undirected binary graph $E'_N$
via
\[
x'_i = \ind_{x_i > 0}\cdot \ind_{i_1 < i_2 < \dots < i_d}, \quad i\in\nats^d_{\neq}.\label{eq:binaryxi}
\]
\cref{cor:undirect} provides the precise statement of the result; note that the bound is identical to that from
\cref{lem:marginal}.
\bncor\label{cor:undirect}
Let $E'_N \defined (x'_i)_{i\in\nats^d_{\neq}} \in\mcN_d$ be the set of edges of a network with truncation $E'_{N,K}\in\mcN_d$,
and denote their marginal distributions $P'_{N}$ and $P'_{N,K}$.
If there exists a measurable function $f$ such that
\[
E'_N = f(E_N)\qquad\text{and}\qquad E'_{N,K} = f(E_{N,K}),
\]
then
\[
\tvd{P'_{N}}{P'_{N,K}} \le 1- \Pr\left( I_N \leq K\right).
\] 
\encor

\subsection{Independent likelihood process}\label{sec:crmtrunc}
We now specialize \cref{lem:marginal} to the setting where $\Theta$ is a CRM
generated by a series representation of the form \cref{eq:seriesrep}, and the network is generated via the independent
likelihood process from \cref{eq:multiedge}.
As a first step towards a bound on the truncation error
for general hypergraphs with $d>1$ in \cref{thm:truncerror}, we present
a simpler corollary in the case where $d = 2$, which
 is of direct interest in analyzing the truncation error
of popular edge-exchangeable networks. 

\bncor\label{cor:1}
In the setting of \cref{lem:marginal}, 
suppose $\Theta$ is a CRM generated from the series representation \cref{eq:seriesrep}, 
edges are generated from  the independent likelihood process \cref{eq:multiedge},
and $d = 2 \leq K$. Then
\[
\Pr\left(I_N \le K\right) \geq \exp\left(- N\cdot B_K\right),
\]
where  
\[
B_{K} &= B_{K,1} + B_{K,2} + B_{K,3}\label{eq:bnk2} \\
B_{K,1} &=  \EE\left[\int_{\reals_+^2}-\log\pi\left(\tau(U_1, \gamma_1+\Gamma_K)\tau(U_2, \gamma_2+\Gamma_K)\right)\dee \gamma_1\dee \gamma_2 \right]\\
B_{K,2} &= 2\EE\left[\int_{0}^{\Gamma_K}\frac{(K-1)}{\Gamma_K}\int_{\Gamma_K}^{\infty}-\log\pi\left(\tau(U_1, \gamma_1)\tau(U_2, \gamma_2)\right)\dee \gamma_1\dee \gamma_2\right]\\
B_{K,3} &= 2\EE\left[\int_{\reals_+}-\log\pi\left(\tau(U_1, \Gamma_K) \tau(U_2, \gamma +\Gamma_K)\right)\dee \gamma\right].
\]
\encor

The proof of this result (and \cref{thm:truncerror} below) in \cref{sec:truncproofs} 
follows by representing the 
tail of the CRM as a unit-rate Poisson process on $[\Gamma_K, \infty)$. 
Though perhaps complicated at first glance, an intuitive interpretation of the 
truncation error terms $B_{K,i}$ 
is provided by \cref{fig:trunc2}.  $B_{K,1}$ corresponds to the
truncation error arising from the upper right quadrant, where both vertices participating in the edge
were in the discarded tail region.
$B_{K,2}$ is the truncation error arising from the bottom right
and upper left quadrants, where one of the two vertices participating in the edge was in the truncation, and the
other was in the tail. Finally, $B_{K, 3}$ represents the truncation
error arising from edges in which one vertex was at the boundary of tail and truncation, and the other was in the tail.


\cref{thm:truncerror} is the generalization of \cref{cor:1} from $d=2$ to the general setting of arbitrary hypergraphs with $d > 1$. 
The bound is analogous to that in \cref{cor:1}---with $B_{K}$ expressed as a sum of terms, each
corresponding to whether vertices were in the tail, boundary, or truncation region---except that
there are $d > 1$ vertices participating in each edge, resulting in more terms in the sum.
Note that \cref{thm:truncerror} also guarantees that the bound is not vacuous, and indeed 
converges to 0 as the truncation level $K\rightarrow \infty$ as expected. 

\bnthm\label{thm:truncerror}
In the setting of \cref{lem:marginal}, 
suppose $\Theta$ is a CRM generated from the series representation \cref{eq:seriesrep}, 
edges are generated from the independent likelihood process \cref{eq:multiedge},
and $1 < d \leq K$. Then
\[
\Pr\left(I_N \le K \right) \geq \exp\left(-N \cdot B_K\right),
\]
 where
\[
B_{K} &= B_{K,1} + B_{K,2} + B_{K,3},
\label{eq:bnkd}\\
B_{K,1} &=  \EE\left[\int_{[\Gamma_K, \infty)^{d}}-\log\pi(\ttheta)\dee\gamma\right]\\
B_{K,2} &= \sum_{\ell=1}^{d-1}\binom{d}{\ell}\EE\!\left[\frac{(K-1)\,!}{(K-1-\ell)\,!} \Gamma_K^{-\ell}\!\!\int_{[0,\Gamma_K]^\ell\times[\Gamma_K,\infty)^{d-\ell}}\!\!-\log\pi(\ttheta)\dee\gamma\right]\\
B_{K,3} &= \sum_{\ell=1}^{d-1}\ell\binom{d}{\ell}\EE\!\left[\frac{(K-1)\,!}{(K-\ell)\,!} \Gamma_K^{-(\ell-1)}\!\!\int_{[0,\Gamma_K]^\ell\times[\Gamma_K,\infty)^{d-\ell}}\!\!-\delta_{\gamma_\ell = \Gamma_K}\log\pi(\ttheta)\dee\gamma\right],
\]
$\delta_\cdot$ is the Dirac delta,
 $\dee \gamma\defined \prod_{j=1}^d \dee\gamma_j$, and $\tilde{\theta} \defined \prod_{j=1}^{d}\tau(U_j, \gamma_j)$.
 Furthermore, $\emph{lim}_{K\rightarrow\infty}B_{K} = 0$.
\enthm
The same geometric intuition from the $d=2$-dimensional case applies
to the general hypergraph truncation error in \cref{eq:bnkd}.
$B_{K,1}$ corresponds to the error arising from the edges whose vertices all belong to
the tail region $\Theta_{K+}$. Each term in the summation in $B_{K,2}$
corresponds to the error arising from edges that have $\ell$ out of $d$ vertices
belonging to the truncation $\Theta_K$. Each term in the summation in
$B_{K,3}$ corresponds to the error arising from the edges that have $\ell-1$ out
of $d$ vertices belonging to the truncation $\Theta_K$ and have one vertex
exactly on the boundary of the truncation. Note that we obtain \cref{cor:1} 
by taking $d=2$ in \cref{eq:bnkd}. 

\subsection{Categorical likelihood process}\label{sec:ncrmtrunc}
We may also specialize \cref{lem:marginal} to the setting where the network is
generated via the single-edge-per-step categorical likelihood process in \cref{eq:categorical}.
However, truncation with the categorical likelihood process poses a few key challenges.
From a practical angle, certain choices of series representation for generating $\Theta$
may be problematic. For instance, when using the rejection representation \cref{eq:rejecrep}
in the typical case where $\mu \neq \nu$, there is a nonzero probability that
\[
\sum_{k=1}^K \ind\left[\theta_k > 0\right] < d,
\]
meaning there aren't enough accepted vertices in the truncation to generate a single
edge. In this case, the categorical likelihood process---which must generate exactly one edge per step---is ill-defined.
An additional theoretical challenge arises from the normalization of the original and truncated networks in \cref{eq:categorical},
which prevents the use of the usual theoretical tools for analyzing CRMs.

Fortunately, the inverse L\'evy representation provides an avenue to address both issues.
The rates $\theta_k$ are all guaranteed to be nonzero---meaning as long as $K \geq d$, the categorical
likelihood process is well-defined---and are decreasing, which enables our theoretical analysis in \cref{sec:prftruncnorm}. However, as mentioned earlier, the inverse L\'evy representation is well-known
to be intractable to use in most practical settings. 

\cref{thm:levy} provides a solution: we use the rejection representation to simulate
the rates $\theta_k$, but instead of simulating for iterations $k=1, \dots, K$, we simulate
\emph{until we obtain $K$ nonzero rates}. This is no longer a sample of a truncated rejection
representation; but \cref{thm:levy} shows that the first $K$ nonzero rates have the same
distribution as simulating $K$ iterations of the inverse L\'evy representation. 
Therefore, we can tailor the analysis of truncation error for the categorical likelihood process
in \cref{thm:truncnorm} to the inverse L\'evy representation, and simulate its truncation for any $K$ using the tractable rejection representation in practice.
\bnthm\label{thm:levy}
Let $\theta_1, \dots, \theta_K$ be the first $K$ rates from the inverse L\'evy representation of a CRM, and let
$\theta'_1, \dots, \theta'_K$ be the first $K$ nonzero rates from any rejection representation of the same CRM. Then
\[
(\theta_1, \dots, \theta_K) \eqd (\theta'_1, \dots, \theta'_K).
\]
\enthm

\bnthm\label{thm:truncnorm} 
In the setting of \cref{lem:marginal}, 
suppose $\Theta$ is a CRM generated from the inverse L\'evy representation \cref{eq:levyrep}, 
edges are generated from the categorical likelihood process \cref{eq:categorical},
and $1 < d \leq K$. Then
\[
\Pr\left(I_N \leq K\right) \geq (1-B_K)^{Nd} \geq 0,
\label{eq:bkd}
\]
where 
\[
B_K \defined \EE\left[\int_{-\infty}^{\infty}Q(\Gamma_K, x)\left(\int_{0}^{1}Q(\Gamma_K u, x)\dee u\right)^{K-d}\left(\frac{\dee}{\dee x}e^{\int_{0}^{\infty}Q(\Gamma_K+\gamma,\, x) - 1\, \dee\gamma}\right) \dee x\right],
\label{eq:bkncrm}
\]
and
\[
Q(u, t) = \exp\left(-\nuinv(u)e^{-t}\right) \quad\text{and} \quad\Gamma_K\sim \distGam(K, 1).
\]
Furthermore, $\lim_{K\rightarrow\infty}B_K = 0$.
\enthm

\subsection{Examples}\label{sec:examples}
We now apply the results of this section to binary undirected
networks ($d=2$) constructed via \cref{eq:binaryxi} from common 
edge-exchangeable networks \cite{cai2016edge}. 
We derive the convergence rate of truncation error, and provide simulations
of the expected number of edges and vertices under the infinite and truncated network.
In each simulation we use
the rejection representation to simulate $\Theta$,
and run
$N=10,000$ steps of network construction.

\begin{figure}[t!]
\begin{subfigure}{0.45\textwidth}
    \centering\includegraphics[width=0.8\textwidth]{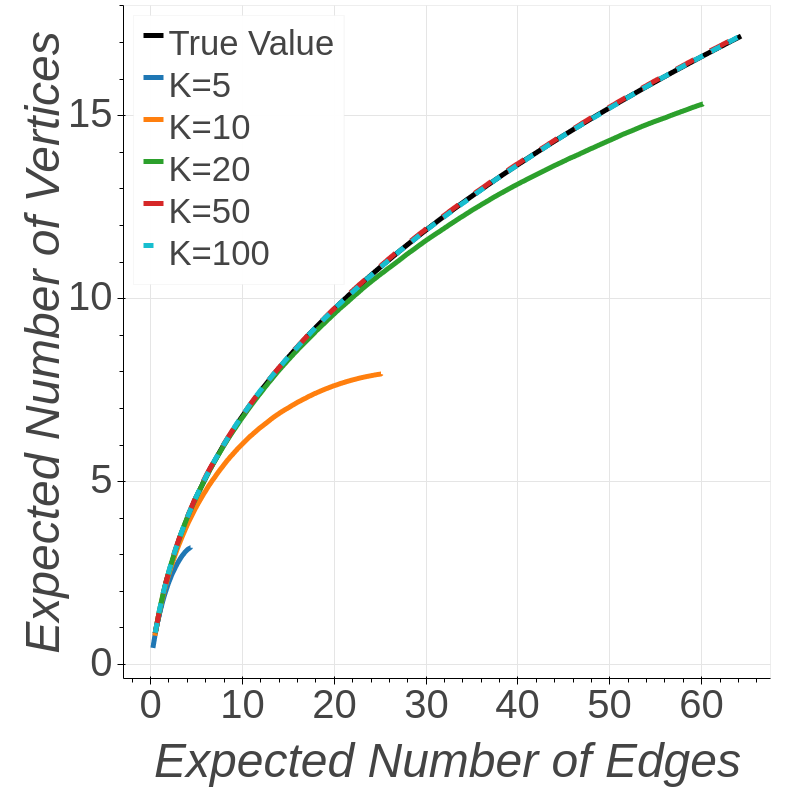}
\caption{Dense graph with $\alpha=0$}\label{fig:beta0}
\end{subfigure}
\begin{subfigure}{0.45\textwidth}
    \centering\includegraphics[width=0.8\textwidth]{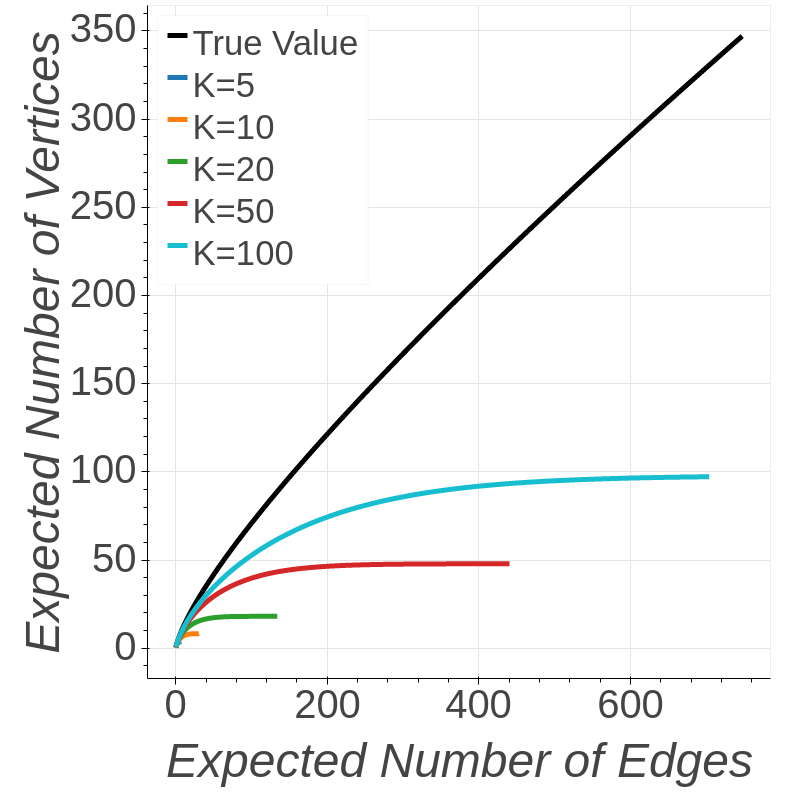}
\caption{Sparse graph with $\alpha=0.6$}\label{fig:betad}
\end{subfigure}
\caption{Beta-independent Bernoulli network}\label{fig:beta}
\end{figure}

\paragraph{Beta-Bernoulli process network}
Suppose $\Theta$ is generated from a beta process,
and network $E_N$ is generated using the independent Bernoulli likelihood process \cite{cai2016edge}. 
The beta process $\distBP(\gamma, \lambda, \alpha)$ \cite{hjort1990} with
discount parameter $\alpha\in[0,1)$, concentration parameter $\lambda>-\alpha$,
and mass parameter $\gamma>0$ is a CRM with rate measure 
\[
\nu(\dee\theta) = \gamma\frac{\Gamma(\lambda + 1)}{\Gamma(1-\alpha)\Gamma(\lambda+\alpha)}\mathds{1}[\theta\le 1]\theta^{-1-\alpha}(1-\theta)^{\lambda+\alpha-1}\dee\theta.
\label{eq:beta}
\] 
The Bernoulli likelihood has the form
\[
 h(x|\theta) = \mathds{1}[x\le1]\theta^{x}(1-\theta)^{1-x}.
\label{eq:bernoulli}
\]
To simulate the process, we use a proposal rate measure $\mu$ given by
\[
\mu(\mathrm{d}\theta) =
\gamma'\ind\left[\theta\leq 1\right]\theta^{-1-\alpha}\dee\theta, \quad \gamma' =
\gamma\frac{\Gamma(\lambda+1)}{\Gamma(1-\alpha)\Gamma(\lambda+\alpha)}.
\]
\emph{Dense network:}
When $\alpha = 0$, the binary beta-Bernoulli graph is dense and 
\[
\muinv(u) = e^{-u/\gamma'}, & & \frac{\dee\nu}{\dee\mu} = (1-\theta)^{\lambda - 1}.
\]
Therefore the rejection representation \cref{eq:rejecrep} of $\distBP(\gamma, \lambda, 0)$ can be written as 
\[
\theta_k = T_k\mathds{1}\left(U_k\le(1-T_k)^{\lambda-1} \right), \quad T_k = e^{-\Gamma_k/\gamma'}.
\]
In \cref{pf:beta}, we show that there exists $K_0\in \nats$ such that
\[
 \forall K \ge K_0, \quad B_{K} \le 12\gamma(1-e^{-1})^{\lambda-2}\left(\frac{\gamma'}{1+\gamma'} \right)^{K}.
\]
This implies that the truncation error of the dense binary beta-Bernoulli network 
converges to 0 geometrically in $K$. \cref{fig:beta0} corroborates this result
in simulation with $\lambda=2$ and $\gamma=1$; 
it can be seen that for the dense beta-Bernoulli graph, truncated graphs with relatively
low truncation level---in this case, $K\approx 50$---approximate the true network model well.

\emph{Sparse network:}
When $\alpha \in (0, 1)$, the binary beta-Bernoulli graph is sparse and 
\[
\muinv(u) = \left(1+ \frac{\alpha u}{\gamma'} \right)^{-1/\alpha} , & & \frac{\mathrm{d}\nu}{\mathrm{d}\mu} = (1-\theta)^{\lambda+\alpha -1}. 
\]
Therefore the rejection representation \cref{eq:rejecrep} of $\distBP(\gamma, \lambda, \alpha)$ can be written as
\[
\theta_k = T_k\mathds{1}\left(U_k\le (1-T_k)^{\lambda+\alpha -1} \right), & & T_k = (1+\alpha\Gamma_k/\gamma')^{-1/\alpha}. 
\]
In \cref{pf:beta}, we show that there exists $K_0\in \nats$ such that
\[
 \forall K \ge K_0, \quad B_{K}\le 6\alpha(\gamma'\alpha^{-1})^{1/\alpha}\  e^{\gamma'\alpha^{-1}}\ (K-1)^{\frac{\alpha-1}{\alpha}}.
\]
This bound suggests that the truncation error 
for the sparse binary beta-Bernoulli network converges to 0
much more slowly than for the dense graph. 
\cref{fig:betad} corroborates this result in simulation with $\lambda = 2$,
$\gamma = 1$, and $\alpha=0.6$; it can be seen that for the sparse beta-Bernoulli
graph, truncated graphs behave significantly differently from the true graph
for moderate truncation levels. In practice, one should select an appropriate (large) 
value of $K$ using our error bounds as guidance, and use sparse data structures to avoid undue memory burden.

\paragraph{Gamma-independent Poisson network}\label{sec:multigammaexample}
Next, consider the network with $\Theta$ generated from a gamma process,
and the network $E_N$ generated using the independent Poisson likelihood process. 
The gamma process $\distGammaP(\gamma, \lambda, \alpha)$
\cite{brix1999generalized} with discount parameter $\alpha\in[0,1)$, scale
parameter $\lambda > 0$ and mass parameter $\gamma>0$ has rate measure 
\[
\nu(\dee\theta) = \gamma\frac{\lambda^{1-\alpha}}{\Gamma(1-\alpha)}\theta^{-\alpha-1}e^{-\lambda\theta}\dee\theta.
\]
The Poisson likelihood has the form
\[
 h(x|\theta) = \frac{\theta^x}{x!}e^{-\theta}.
\]

\begin{figure}[t!]
\begin{subfigure}{0.45\textwidth}
    \centering\includegraphics[width=0.8\textwidth]{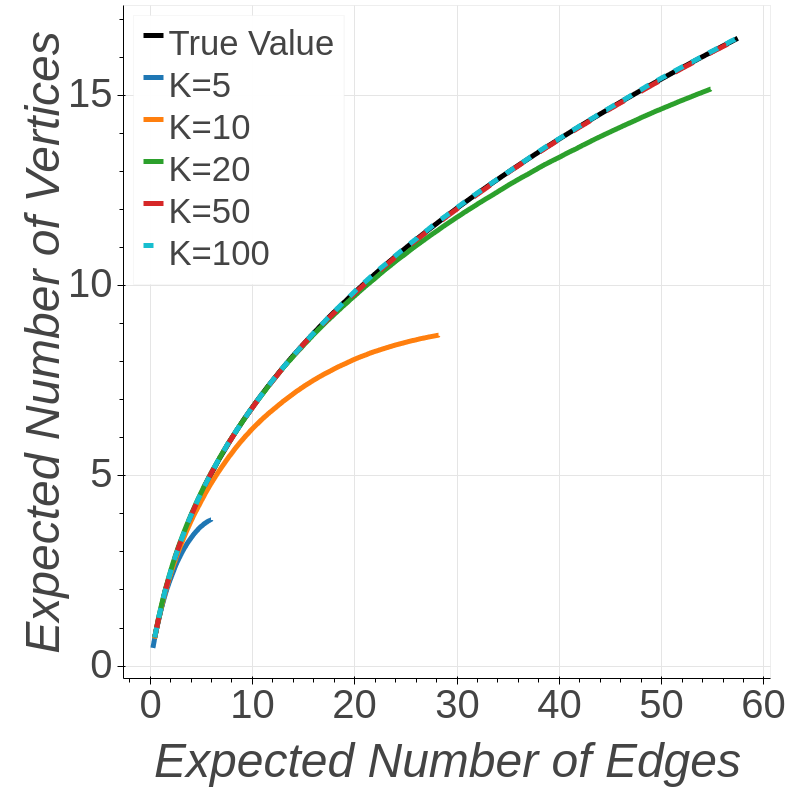}
\caption{Dense graph with $\alpha=0$}\label{fig:gamma0}
\end{subfigure}
\begin{subfigure}{0.45\textwidth}
    \centering\includegraphics[width=0.8\textwidth]{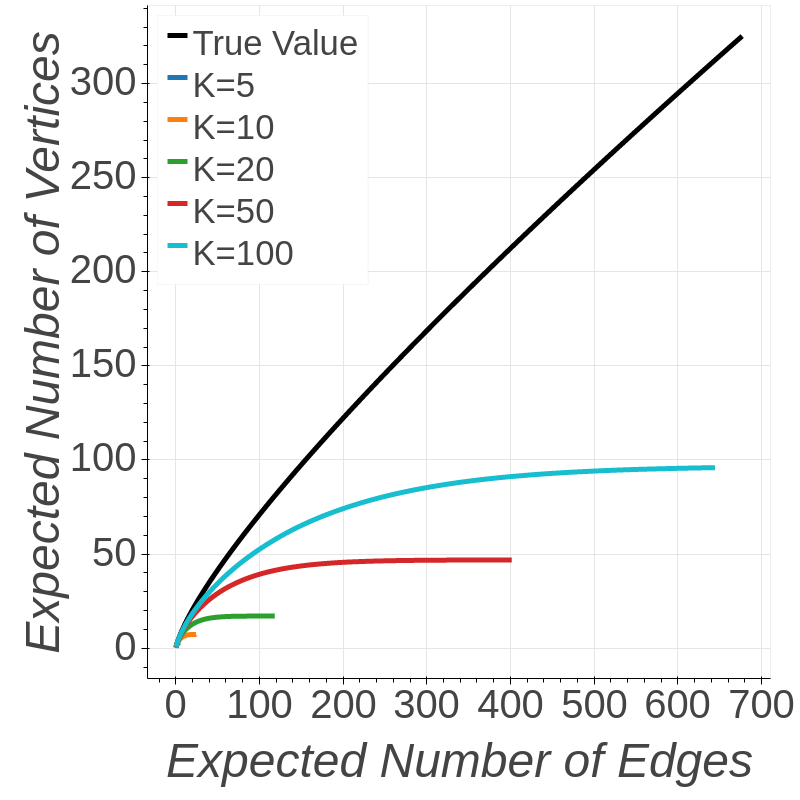}
\caption{Sparse graph with $\alpha=0.6$}\label{fig:gammad}
\end{subfigure}
\caption{Gamma-independent Poisson graph}\label{fig:gamma}
\end{figure}

\emph{Dense network:}
When $\alpha=0$, the gamma-Poisson graph is dense, and we choose the proposal measure $\mu$ to be
\[
\mu(\dee\theta) = \gamma\lambda\theta^{-1}(1+\lambda\theta)^{-1}\mathrm{d}\theta,
\]
such that
\[
 \muinv(u) = 1/\left(\lambda\left(e^{(\gamma\lambda)^{-1}u} - 1\right)\right), \qquad \frac{\mathrm{d}\nu}{\mathrm{d}\mu} = (1+\lambda\theta)e^{-\lambda\theta}.
\]
 Therefore, the rejection representation in \cref{eq:rejecrep} has the form
\[
\theta_k = T_k\mathds{1}\left(U_k\le(1+\lambda T_k)e^{-\lambda T_k}\right), & &  T_k = \frac{1}{\lambda\left(e^{(\gamma\lambda)^{-1}\Gamma_k} - 1 \right)}.
\]
In \cref{pf:gamma}, we show that there exists $K_0\in \nats$ such that
\[
 \forall K \ge K_0, \quad B_{K} \le 6\frac{\gamma}{\lambda}\left(\frac{\gamma\lambda}{1+\gamma\lambda} \right)^{K-1}.
\]
Again, for the dense network $B_{K}$ converges to 0 geometrically, indicating that truncation may provide reasonable
approximations to the original network. \cref{fig:gamma0} corroborates this result in simulation with
$\lambda = 2$ and $\gamma=1$; for $K\approx 50$, no vertices are discarded on average by truncation. 

\emph{Sparse network:}
When $\alpha \in (0, 1)$, the gamma-Poisson graph is sparse, and we choose
the proposal measure $\mu$ to be
\[
\mu(\dee\theta) = \gamma\frac{\lambda^{1-\alpha}}{\Gamma(1-\alpha)}\theta^{-1-\alpha}\mathrm{d}\theta,
\]
such that
\[
\muinv(u) = (\gamma'u^{-1})^{1/d}, \quad \gamma'\defined \gamma\frac{\lambda^{1-\alpha}}{\alpha\Gamma(1-\alpha)}, \quad
\frac{\dee\nu}{\mathrm{d}\mu} = e^{-\lambda\theta}.
\]
Therefore the rejection representation in \cref{eq:rejecrep} has the form
\[
\theta_k = T_k\mathds{1}\left(U_k\le e^{-\lambda T_k} \right), \quad T_k = \left(\gamma'\Gamma_k^{-1} \right)^{1/\alpha}.
\]
In \cref{pf:gamma}, we show there exists $K_0\in \nats$ such that
\[
 \forall K\ge K_0, \quad B_{K}\le \frac{12\gamma^2\lambda^{1-\alpha}}{(1-\alpha)\Gamma(1-\alpha)}\left(\frac{3\gamma'}{K-1}\right)^{\frac{1-\alpha}{\alpha}}.
\]
Again, for the sparse network $B_K$ converges to 0 slowly, suggesting that the truncation error for the sparse
binary gamma-Poisson graph converges more slowly than for the dense graph. \cref{fig:gammad} corroborates this
result in simulation with $\lambda = 2$, $\gamma = 1$, and $\alpha = 0.6$;
for a moderate range of truncation values $K\leq 100$, the truncated graph behaves very differently from the true graph.
Therefore in practice, one should follow the above guidance for the beta-Bernoulli network:
select a value of $K$ using our error bounds, and avoid intractable memory requirements by using sparse data structures.

\section{Truncated posterior inference}\label{sec:inference}
In this section, we develop a tractable approximate posterior inference method for network
models via truncation, and analyze its approximation error.
In particular, given an observed network $E_N$, we want to simulate from
the posterior distribution of the CRM rates and the parameters of the CRM rate measure. 
Since exact expressions of full conditional densities are not available, 
we truncate the model and run an approximate Markov chain Monte Carlo algorithm.
We provide a rigorous theoretical justification for this simple approach by establishing a bound
on the total variation distance between the truncated and exact posterior distribution. 
This in turn provides a method to select the truncation level in a principled manner.

Although this section focuses on self-product-CRM-based network models, the method
and theory we develop are both general and extend easily to other CRM-based models, e.g.~for 
clustering \cite{Antoniak74}, feature allocation \cite{Griffiths06}, and trait allocation \cite{campbell2018exchangeable}. In 
particular, the methodology only requires 
 bounds on tail occupancy probabilities (e.g., in the present context, the probability that $I_N \leq K$)
rather than the exact evaluation of these quantities.

\subsection{Truncation error bound}
We begin by examining the density of the
posterior distribution 
of the $K^\text{th}$ rate from the inverse L\'evy representation $\theta_K$ for some fixed $K\in\nats$,
the unordered collection of rates $(\theta_k)_{k=1}^{K-1}$ such that $\theta_k \geq \theta_K$,
and the parameters $\sigma\in\reals^m$ of the CRM rate measure $\nu_\sigma$
given the observed set of edges $E_N$.
We denote $\nu_\sigma(\theta)$ to be the density of $\nu_\sigma(\dee \theta)$ and $p(\sigma)$ to be the prior density of $\sigma$, 
both with respect to the Lebesgue measure. Given these definitions the posterior density can be expressed as 
\[
&p(\sigma, \theta_{1:K}, X_{1:N}) \propto p(\sigma)\cdot e^{-\nu_\sigma[\theta_K, \infty)}\cdot\prod_{k=1}^K\ind[\theta_K \leq \theta_k]\nu_\sigma(\theta_k)\\
&\hspace{3cm}\cdot \left(\prod_{n=1}^{N}\prod_{i\in[K]^d_{\neq}}h(x_{ni} \given {\scriptstyle \prod_{j=1}^d} \theta_{i_j})\right)
\cdot p(x_{1:N}^{K+} | \theta_{1:K}, \sigma), \label{eq:postdensity}
\]
where $[K]^d_{\neq}$ is the subset of $\nats^d_{\neq}$ such that $\max_{j\in[d]} i_j \leq K$,
and $x_{1:N}^{K+}$ is shorthand for the set of $(x_{ni})_{n\in [N], i\in\nats_{\neq}^d}$ such that $i \notin [K]^d_{\neq}$.
All the factors on the first row correspond to the prior of $\sigma$ and $(\theta_k)_{k=1}^K$, 
and the first factor on the second row corresponds to the likelihood of the portion of the network
involving only vertices $1\dots K$; these are straightforward to evaluate, though $\nu_\sigma[\theta_K, \infty)$
will occasionally require a 1-dimensional numerical integration.
The last factor corresponds 
to the likelihood of the portion of the network involving vertices 
beyond $K$, and is not generally possible to evaluate exactly.

To handle this term, suppose that $K$ is large enough that $x_{1:N}^{K+} = 0$. 
Then
$p(x_{1:N}^{K+} | \theta_{1:K}, \sigma) = \Pr\left(I_{N} \leq K \given \Gamma_{1:K}, \sigma \right)$, i.e.,
 the probability that no portion of the network involves vertices of index $> K$.
\cref{thm:conditionaltrunc} provides upper and lower bounds on this probability akin to those of 
\cref{thm:truncerror}---indeed, \cref{thm:conditionaltrunc} is an intermediate step in the proof of \cref{thm:truncerror}---that
apply conditionally on the values of $U_{1:K}$, $\Gamma_{1:K}$ rather than marginally.
This theorem also makes the dependence of the bound on the rate measure parameters $\sigma$ notationally explicit
via parametrized series representation components $\tau_\sigma$ and $g_\sigma$ from \cref{eq:seriesrep}.
Finally, though \cref{thm:conditionaltrunc} applies to general series representations,
we require only the specific instantiation for the inverse L\'evy representation in the present context.

\bnthm\label{thm:conditionaltrunc}
The conditional probability $\Pr\left(I_N \le K \given U_{1:K}, \Gamma_{1:K},\sigma\right)$ satisfies 
\[
1 \geq \Pr\left(I_N \le K \given U_{1:K}, \Gamma_{1:K}, \sigma \right) \geq \exp\left(-N \cdot B(U_{1:K},\Gamma_{1:K}, \sigma)\right),
\]
where
\[
B(U_{1:K},\Gamma_{1:K}, \sigma)&=
\sum_{\ell=0}^{d-1}\binom{d}{\ell}
\sum_{\begin{subarray}{c}\mcL\subseteq [K]\\ |\mcL|=\ell\end{subarray}}B(U_{1:K}, \Gamma_{1:K}, \sigma, \mcL)\\
 B(U_{1:K}, \Gamma_{1:K}, \sigma, \mcL)&=\int_{[\Gamma_K, \infty)^{d-\ell}}\!\!\!\!\EE\left[-\log\pi\left(\prod_{j\in\mcL}\theta_{j}\prod_{j=1}^{d-\ell}\tau_\sigma(U'_{j}, \gamma_{j})\right) \given \theta_{1:K}\right]\dee\gamma,
\]
where $\dee\gamma = \prod_{j=1}^{d-\ell}\dee \gamma_j$ and $U'_j \distiid g_\sigma$.
\enthm
Since $U_{1:K}$ is unused in the inverse L\'evy representation, 
in the present context we use the notation $B(\Gamma_{1:K},\sigma)$ for brevity.
The bound in \cref{thm:conditionaltrunc} implies that as long as $K$ is set large enough such that 
both $x_{1:N}^{K+} = 0$ and $B(\Gamma_{1:K}, \sigma) \le \epsilon/N$ for some $\epsilon > 0$
then
\[
1-\epsilon \leq p(x_{1:N}^{K+} | \theta_{1:K}, \sigma) \leq 1.
\]
Therefore as $K$ increases, this term should become $\approx 1$ and have a negligible effect on 
the posterior density. We use this intuition to propose a truncated Metropolis--Hastings
algorithm that sets $K > I_N$, ignores the $p(x_{1:N}^{K+} | \theta_{1:K}, \sigma)$ term in the acceptance ratio,
and fixes $x_{N}^{K+}$ to 0.
\cref{thm:mherror} provides a rigorous analysis of the error involved in using the truncated sampler.
\bnthm\label{thm:mherror}
Fix $K > I_N$.
Let $\Pi$ be the distribution of $\sigma, \theta_{1:K}$ given $E_N$,
and let $\hPi$ be the distribution with density proportional to \cref{eq:postdensity} without 
the $p(x_{1:N}^{K+} | \theta_{1:K}, \sigma)$ term and with $x_{1:N}^{K+}$ fixed to 0.
If for some $\eta\in[0,1)$,
\[
\hPi\left\{B(\Gamma_{1:K}, \sigma) \le \epsilon/N\right\} \geq 1-\eta, 
\]
then 
\[
\tvd{\hPi}{\Pi} &\leq \frac{3(\epsilon+\eta)}{2} - \epsilon\eta. \label{eq:posteriorl1bound}
\]
\enthm
\subsection{Adaptive truncated Metropolis--Hastings}
Crucially, as long as one can obtain samples from the
\emph{truncated} posterior distribution, one can estimate the bound in \cref{eq:posteriorl1bound}
using sample estimates of the tail probability $\hPi\left\{B(\Gamma_{1:K}, \sigma) \le \epsilon/N\right\} \geq 1-\eta$.
Therefore, one can compute the bound \cref{eq:posteriorl1bound}
without needing to evaluate $p(x_{1:N}^{K+} | \theta_{1:K}, \sigma)$ exactly.
This suggests the following adaptive truncation procedure:
\benum
\item Pick a value of $K > I_N$ and desired total variation error $\xi \in (0, 1)$.
\item Obtain samples from the truncated posterior.
\item Minimize the bound in \cref{eq:posteriorl1bound} over $\epsilon \in(0,1)$, using the samples to estimate $\eta = 1 - \hPi\left\{B(\Gamma_{1:K}, \sigma) \le \epsilon/N\right\}$ for each value of $\epsilon$.
\item If the minimum bound exceeds $\xi$, increase $K$ and return to 2.
\item Otherwise, return the samples.
\eenum
In this work, we start by initializing $K$ to $I_N+1$.
In order to decide how much to increase $K$ by in each iteration, we use
the sequential representation to extrapolate
the total variation bound \cref{eq:posteriorl1bound} to larger values of $K$ without actually performing MCMC sampling.
In particular, for each posterior sample, we use its hyperparameters to generate additional rates
from the sequential representation. We then use these extended samples to compute the 
total variation error guarantee as per step 3~above. We continue to generate additional rates (typically doubling the number each time)
until the predicted total variation guarantee is below the desired threshold. Finally, we use linear interpolation
between the last two predicted errors to find the next truncation level $K$ that matches the desired (log) error threshold.
This fixes a new value of $K$; at this point, we return to step 2~above and iterate.

\begin{figure}[t!]
	\begin{subfigure}{0.32\textwidth}
		\centering\includegraphics[width=\textwidth]{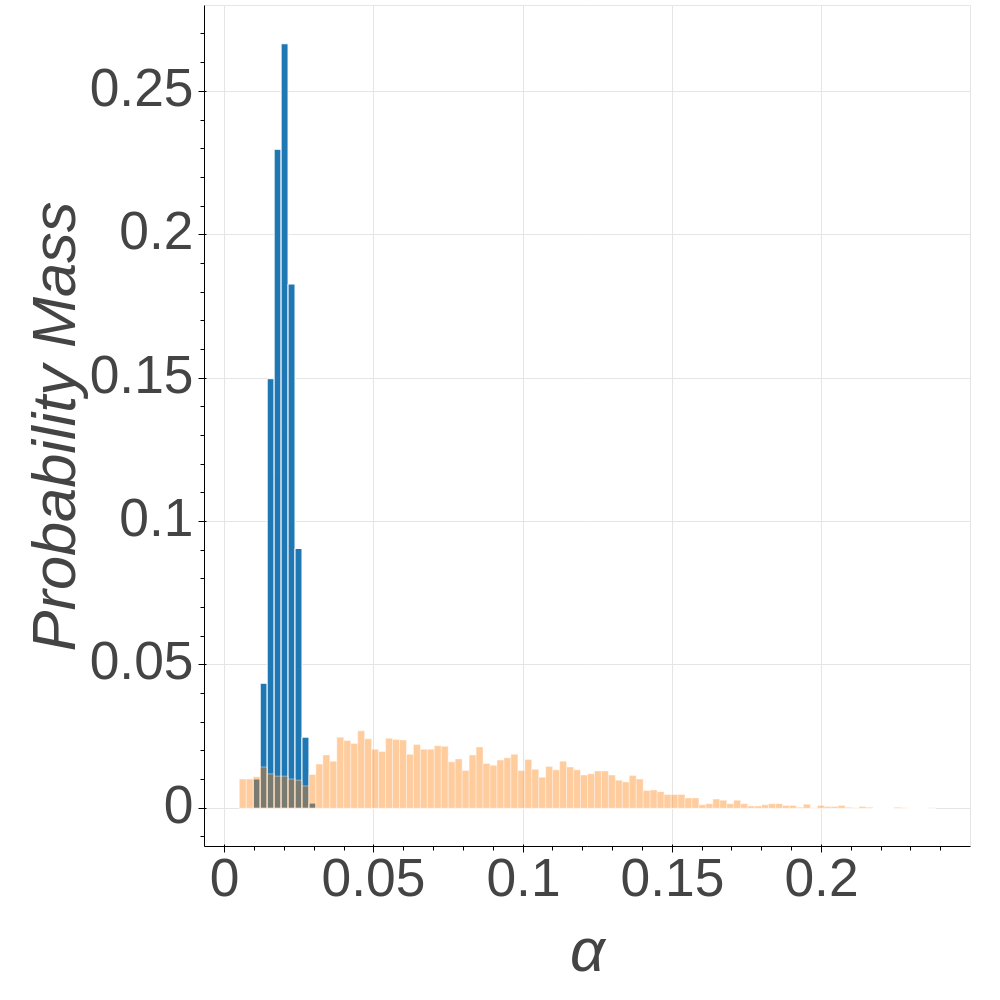}
		\caption{True $\alpha=0$}\label{fig:alph_dense}
	\end{subfigure}
	\begin{subfigure}{0.32\textwidth}
		\centering\includegraphics[width=\textwidth]{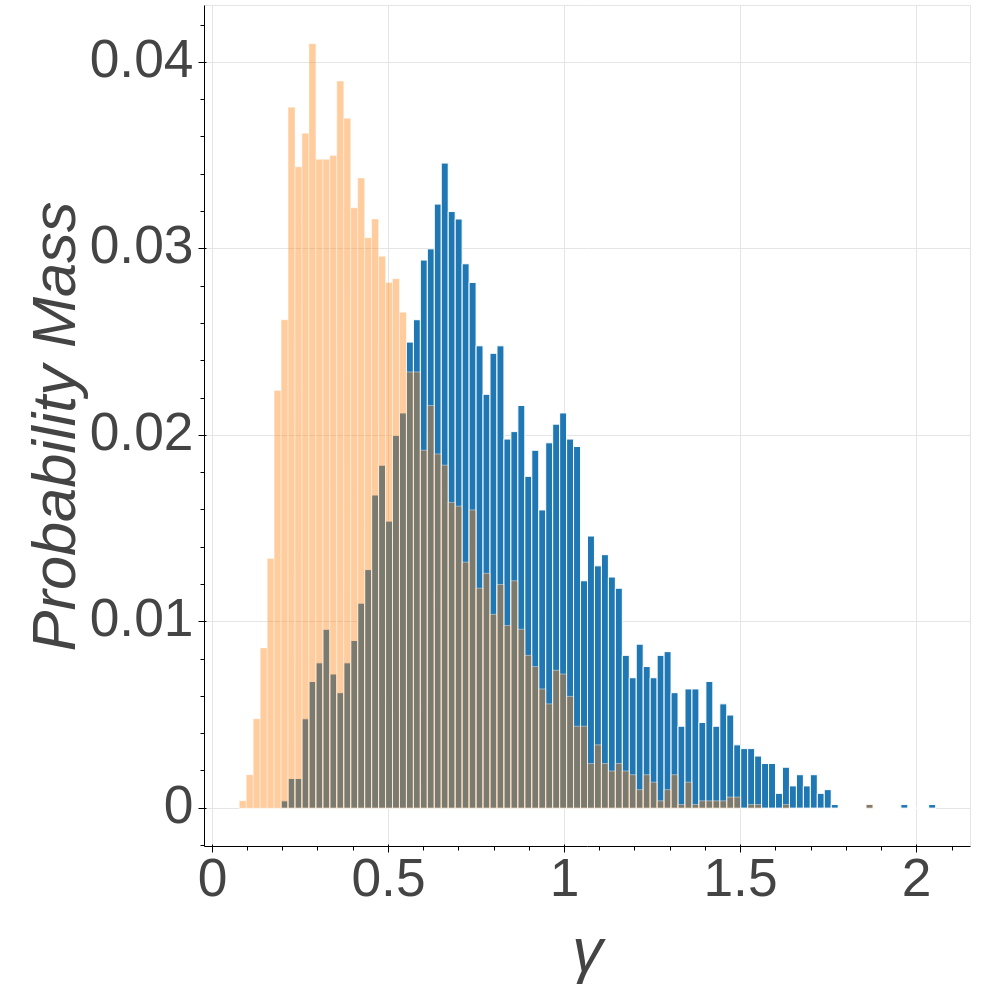}
		\caption{True $\gamma=1$}\label{fig:gam_dense}
	\end{subfigure}
	\begin{subfigure}{0.32\textwidth}
		\centering\includegraphics[width=\textwidth]{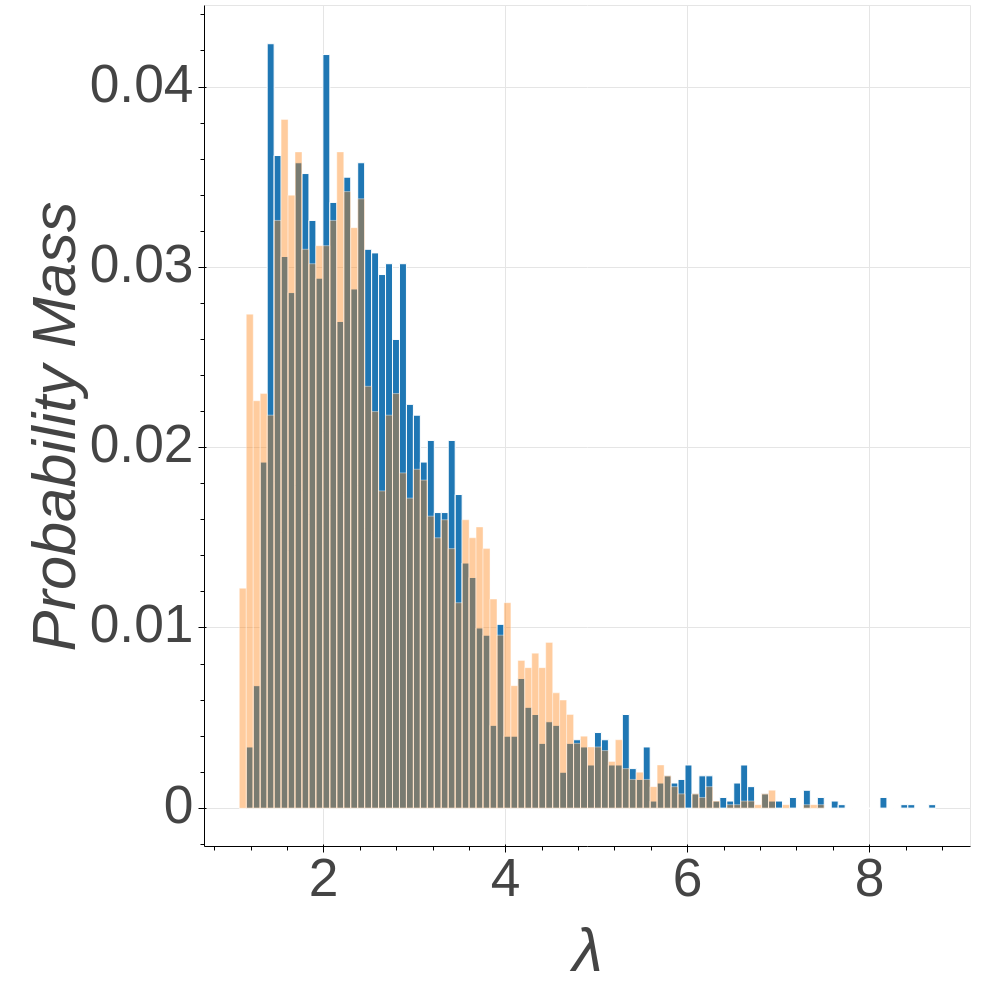}
		\caption{True $\lambda=2$}\label{fig:lamb_dense}
	\end{subfigure}
	\caption{Posterior of $\alpha$, $\gamma$, $\lambda$ for the dense simulated network. Orange histograms
depict the first adaptation iteration, and blue histograms depict the second and final iteration.}\label{fig:syn_dense}
\end{figure}
\begin{figure}[t!]
	\begin{subfigure}{0.32\textwidth}
		\centering\includegraphics[width=\textwidth]{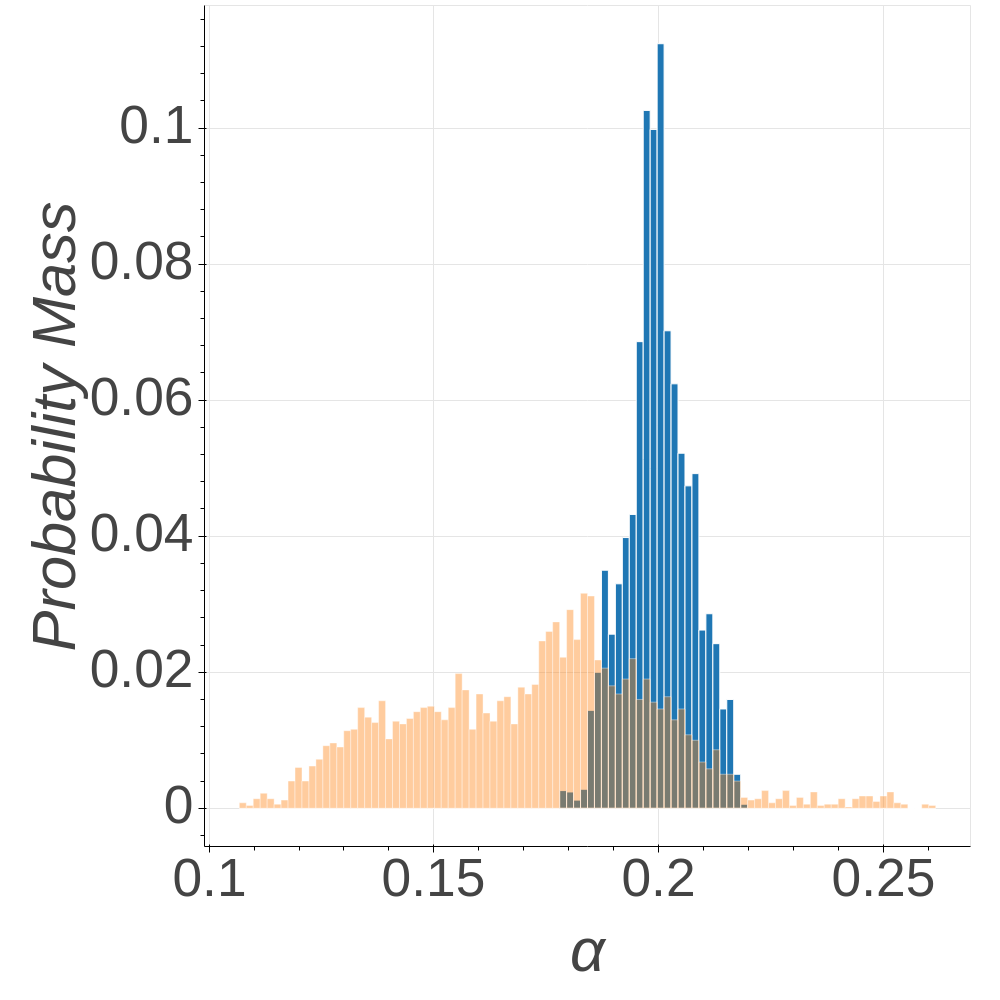}
		\caption{True $\alpha=0.2$}\label{fig:alph_sparse}
	\end{subfigure}
	\begin{subfigure}{0.32\textwidth}
		\centering\includegraphics[width=\textwidth]{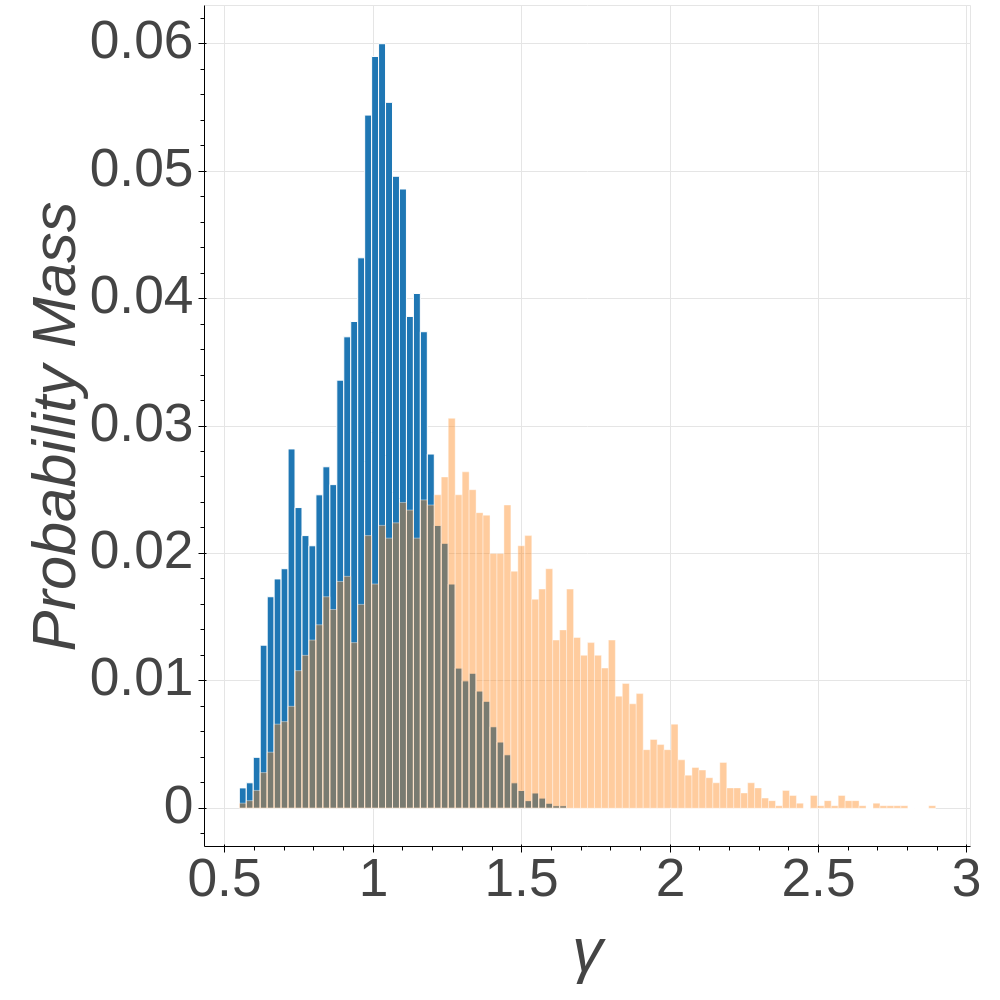}
		\caption{True $\gamma=1$}\label{fig:gam_sparse}
	\end{subfigure}
	\begin{subfigure}{0.32\textwidth}
		\centering\includegraphics[width=\textwidth]{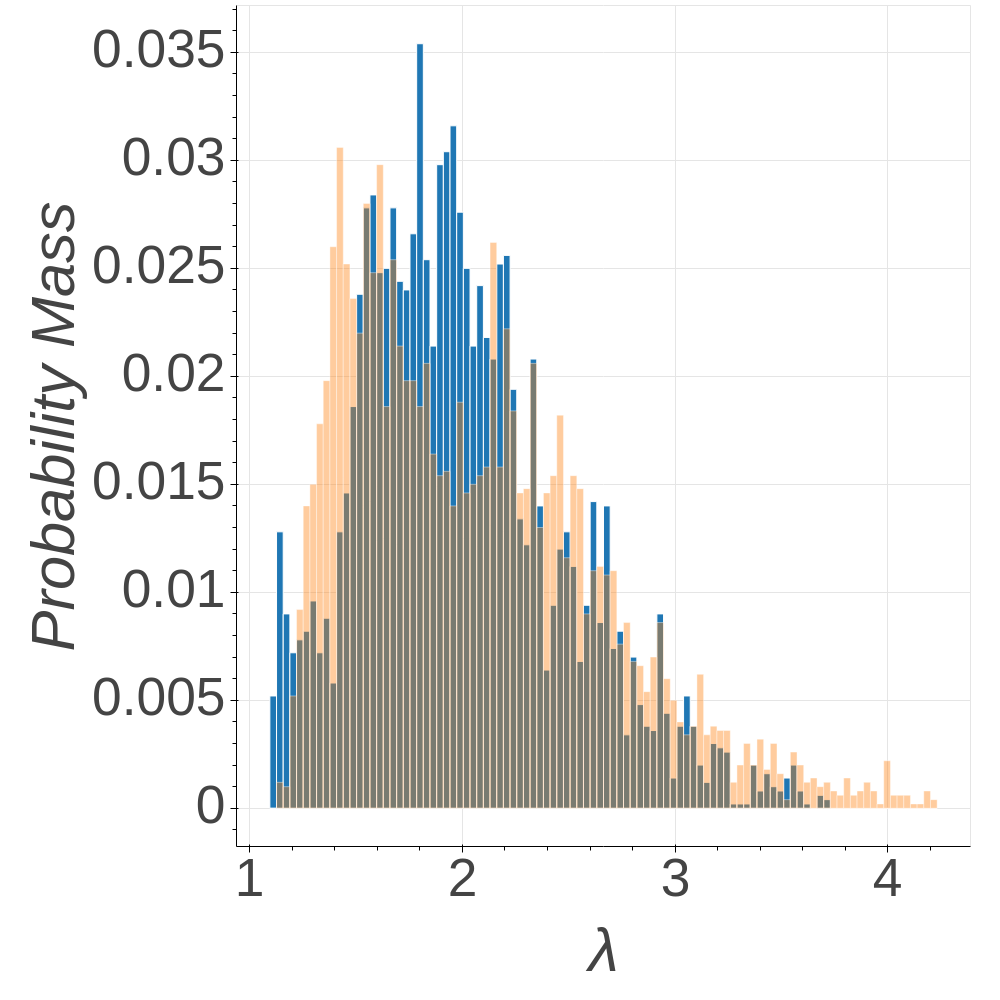}
		\caption{True $\lambda=2$}\label{fig:lamb_sparse}
	\end{subfigure}
	\caption{Posterior of $\alpha$, $\gamma$, $\lambda$ for the sparse simulated network. Orange histograms
depict the first adaptation iteration, and blue histograms depict the second and final iteration.}\label{fig:syn_sparse}
\end{figure}

\subsection{Experiments}
In this section, we examine the properties of the proposed adaptive truncated inference
algorithm for the beta-independent Bernoulli network model in \cref{eq:beta,eq:bernoulli} with discount $\alpha\in(0,1)$,
concentration $\lambda > 1$, mass $\gamma > 0$, unordered collection of rates $\theta_{1:K-1}$, 
and $K^\text{th}$ rate from the sequential representation $\theta_K$. 
In order to simplify inference, we transform each of these parameters to an unconstrained version:
\[
\alpha &= \frac{\exp(\alpha_u)}{1+\exp(\alpha_u)}, &
\lambda &= 1+\exp(\lambda_u), &
\gamma &= \exp(\gamma_u),\\
\theta_K &= \frac{\exp(\theta_{u,K})}{1+\exp(\theta_{u,K})}, &
\theta_k &= \frac{\theta_K + \exp(\theta_{u,k})}{1+\exp(\theta_{u,k})} & k&=1, \dots, K-1.
\]
We use a Markov chain Monte Carlo algorithm that includes 
an exact Gibbs sampling move for $\gamma$,
and separate Gaussian random-walk Metropolis--Hastings moves 
for $\alpha_u$, $\lambda_u$, $\theta_{u,K}$,
all $\theta_{u,k}$ such that vertex $k$ has degree 0 (jointly),
and each $\theta_{u,k}$ such that vertex $k$ has nonzero degree (individually).

\paragraph{Synthetic data}

\begin{figure}[t!]
	\begin{subfigure}{0.48\textwidth}
		\centering\includegraphics[width=\textwidth]{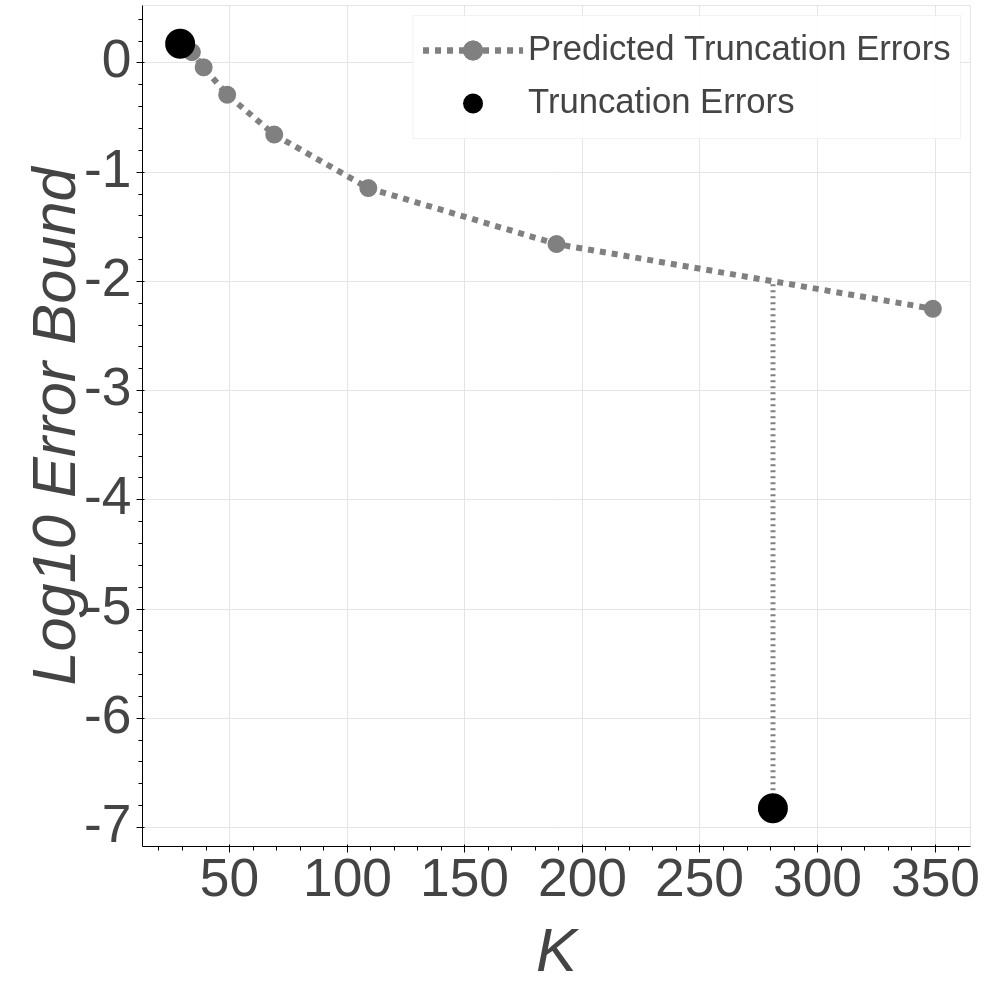}
		\caption{Dense synthetic network}\label{fig:expansion_dense}
	\end{subfigure}
	\begin{subfigure}{0.48\textwidth}
		\centering\includegraphics[width=\textwidth]{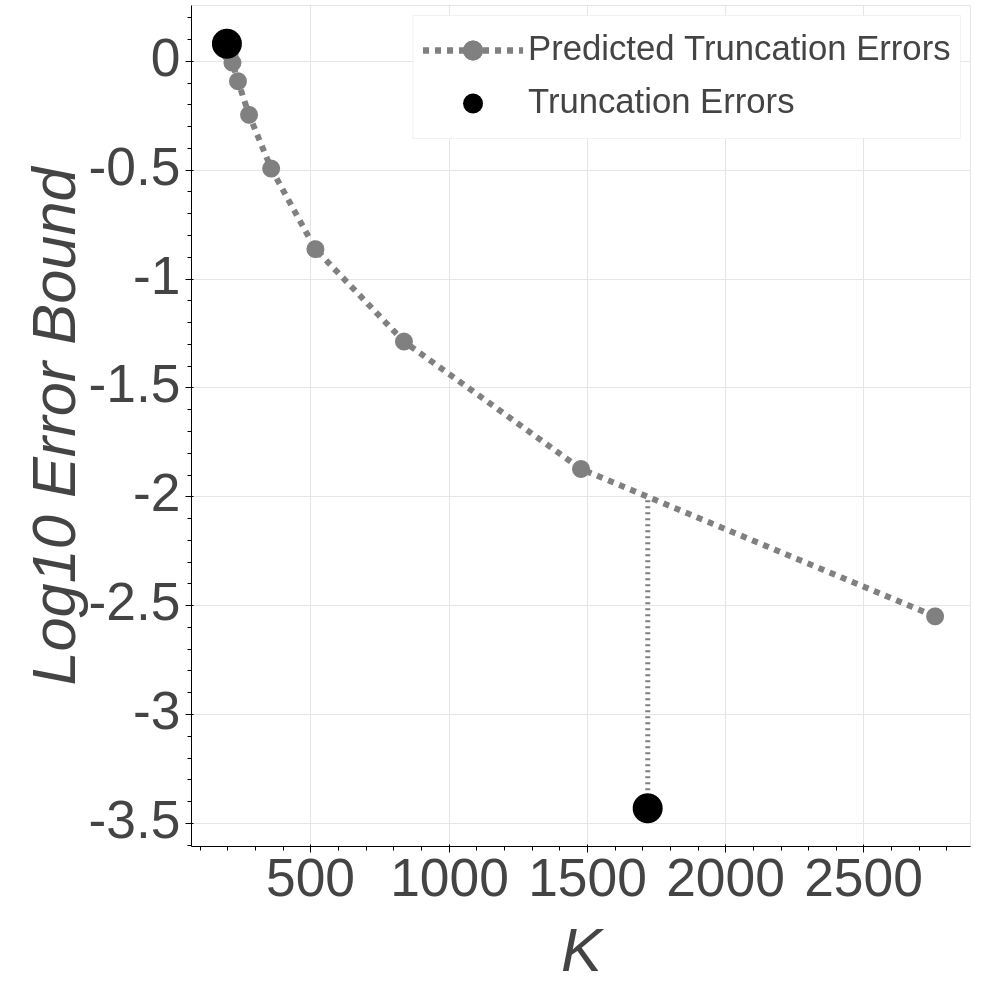}
		\caption{Sparse synthetic network}\label{fig:expansion_sparse}
	\end{subfigure}
	\caption{Visualization of the truncation expansion procedure. Black
circles denote the truncation level and error in each iteration of the
adaptation. The adaptation stops when the truncation error falls below the desired
threshold (here $\log_{10}(0.01) = -2$). Grey dashed lines and circles visualize the predicted truncation
error using rates generated from the sequential representation. The vertical
dotted line shows that the algorithm selects a value of $K$ that attempts to
match the desired $\log_{10}$ error threshold of $-2$ using the
predictions.}\label{fig:expansion}
\end{figure}

We first apply the model to synthetic data simulated from the generative model.
We simulate a sparse network with parameters $\lambda = 2$, $\gamma = 1$, $\alpha=0.2$, 
and $N= 10^5$,
and a dense network with $\lambda = 2$, $\gamma = 1$, $\alpha = 0$, and $N=10^7$. 
In both settings
we set the truncation level for data generation to $500$,
the desired total variation bound from \cref{eq:posteriorl1bound} to $0.01$,
and
initialize the sampler with $\alpha = 0.4$, $\lambda = 5$, $\gamma = 2$ and $\theta$ generated from the sequential representation.
All Metropolis--Hastings moves have proposal standard deviation $0.1$
except the sparse network $\alpha_u$ move, which has standard deivation $0.03$.

\cref{fig:syn_dense,fig:syn_sparse} show histograms of $5{,}000$ marginal
posterior samples of the hyperparameters for the dense (true $\alpha=0$) and
sparse (true $\alpha=0.2$) networks.  In both cases, the approximate posterior
in the first round of adaptation (orange histogram) does not concentrate on the
true hyperparameter values, despite the relatively large number of generative
rounds $N$.  \cref{fig:expansion}---which displays the truncation error and
predictive adaptation procedure---shows why this is the case.  In both
networks, the first adaptation iteration identifies a large truncation error.
After a single round of adaptation, the approximate posterior distributions
(blue histograms) in \cref{fig:syn_dense,fig:syn_sparse} concentrate more on
the true values as expected, and the truncation errors fall well below the
desired threshold ($\log_{10}(0.01) = -2$). It is worth noting that the
predictive extrapolation appears to be quite conservative in these examples,
and especially in the dense network. This happens because the approximate
posterior for the dense network (respectively, sparse
network) assigns mass to higher values of $\alpha$ (respectively, $\gamma$)
than it should, which results in larger truncation error and thus a larger
predicted required truncation level.

\paragraph{Real network data}
Next, we apply the model to a Facebook-like
Social Network\footnote{Available at \url{https://toreopsahl.com/datasets/}}
\cite{panzarasa2009patterns}. The original source network contains 
a sequence of $61{,}734$ weighted,
time-stamped edges, and $1{,}899$ vertices. We preprocess
the data by removing the edge weights, binning the edge sequence into
30-minute intervals, and removing the initial transient of network growth,
resulting in $1{,}899$ vertices and $10{,}435$ edges over
$N=6{,}427$ rounds of generation.
We again set a desired total variation error guarantee of 0.01 during inference.
All Metropolis--Hastings moves have proposal standard deviation $0.1$
except the $\alpha_u$ and degree-0 $\theta_u$ moves, which have standard deviation $0.04$
and $0.03$, respectively. We initialize the sampler to $\alpha=0.1$, $\gamma = 2$, $\lambda = 20$
and sample rates $\theta$ from the prior sequential representation.

\begin{figure}[t!]
	\begin{subfigure}{0.32\textwidth}
		\centering\includegraphics[width=\textwidth]{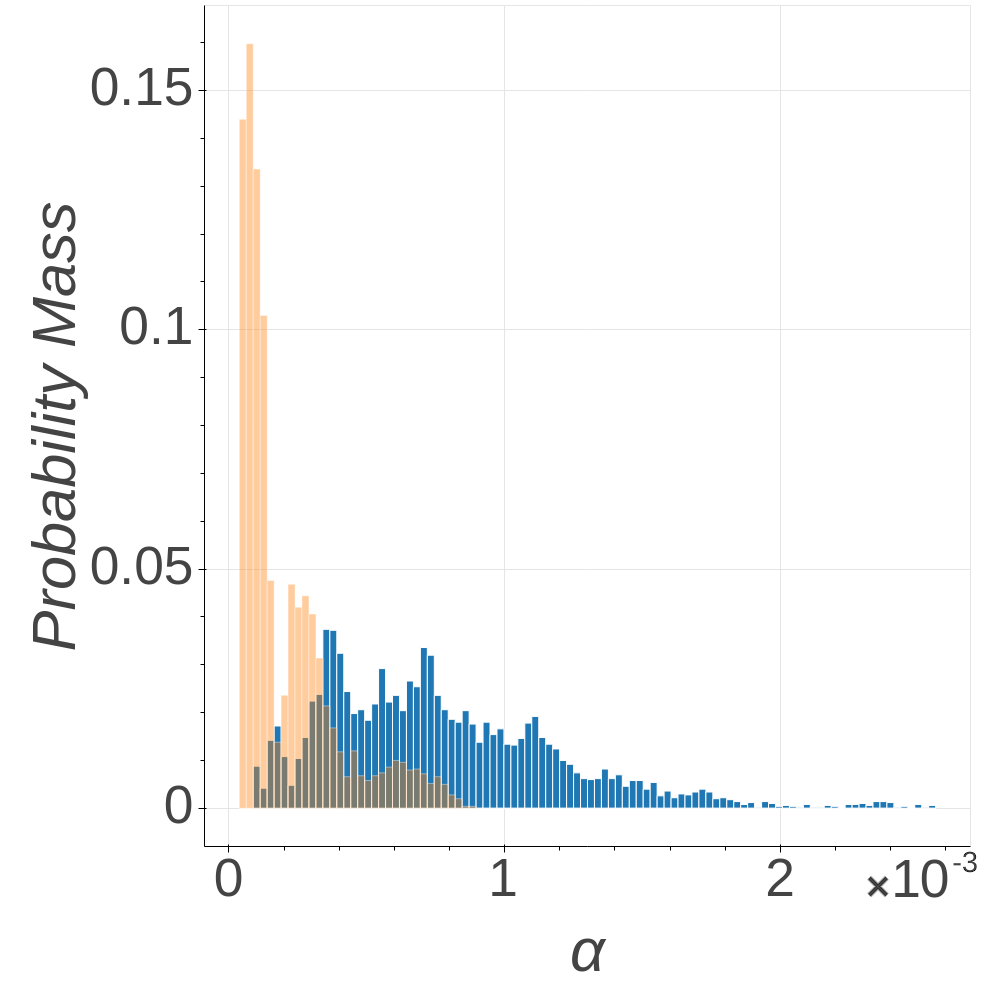}
	\end{subfigure}
	\begin{subfigure}{0.32\textwidth}
		\centering\includegraphics[width=\textwidth]{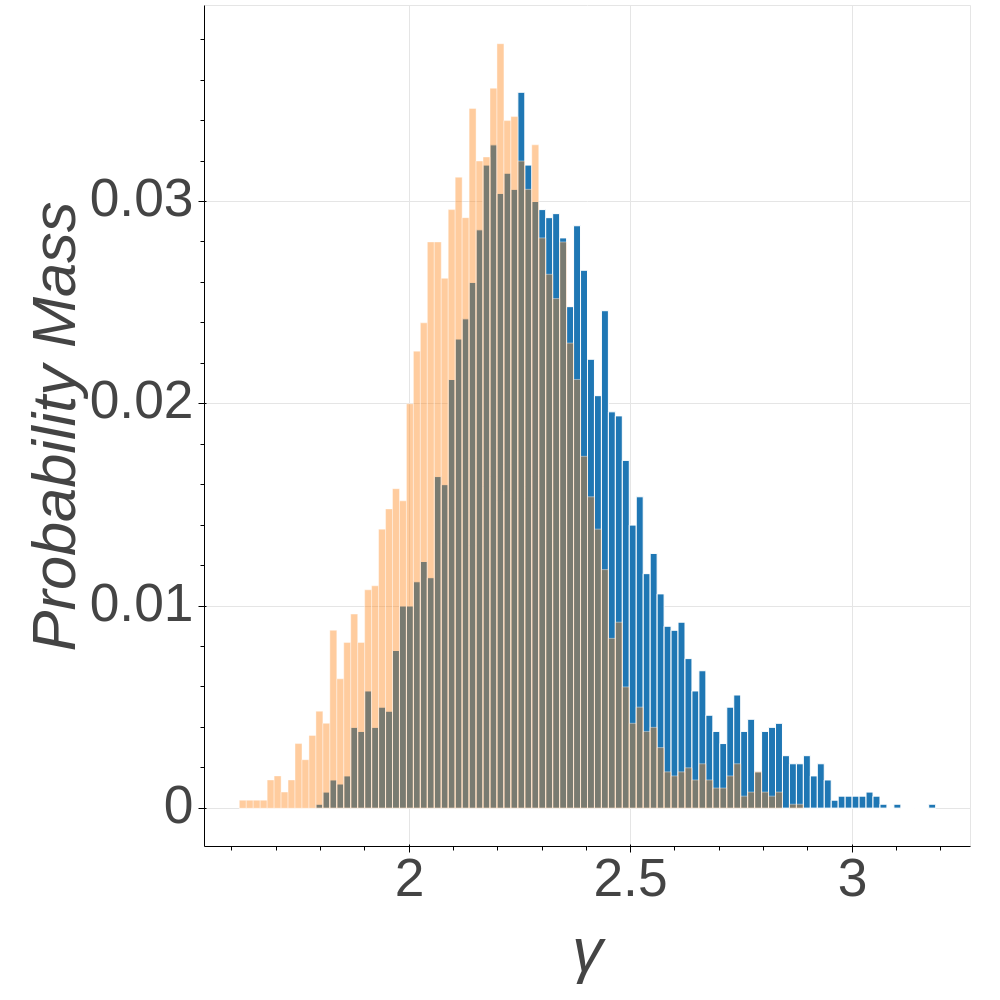}
	\end{subfigure}
	\begin{subfigure}{0.32\textwidth}
		\centering\includegraphics[width=\textwidth]{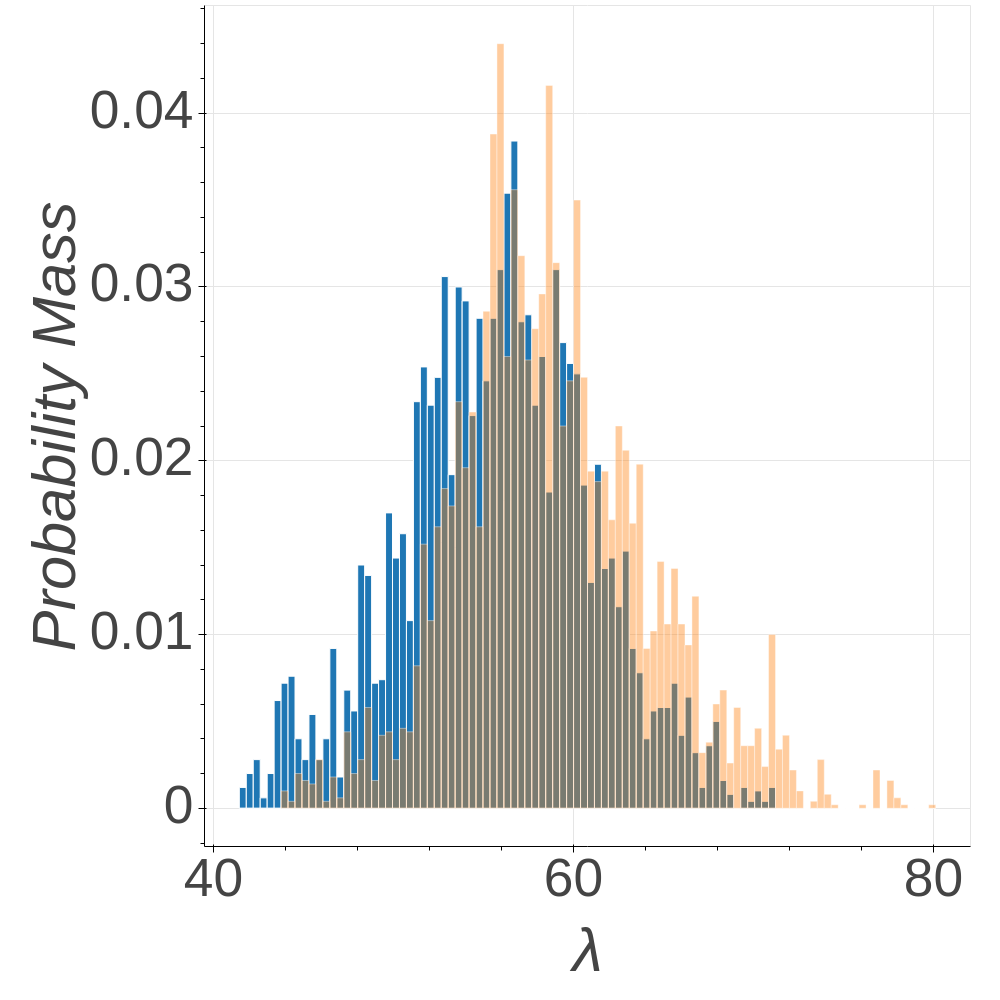}
	\end{subfigure}
	\caption{Posterior of $\alpha$, $\gamma$, $\lambda$ in the Facebook-like network}\label{fig:fb}
\end{figure}

\begin{figure}[t!]
	\begin{subfigure}{0.32\textwidth}
		\centering\includegraphics[width=\textwidth]{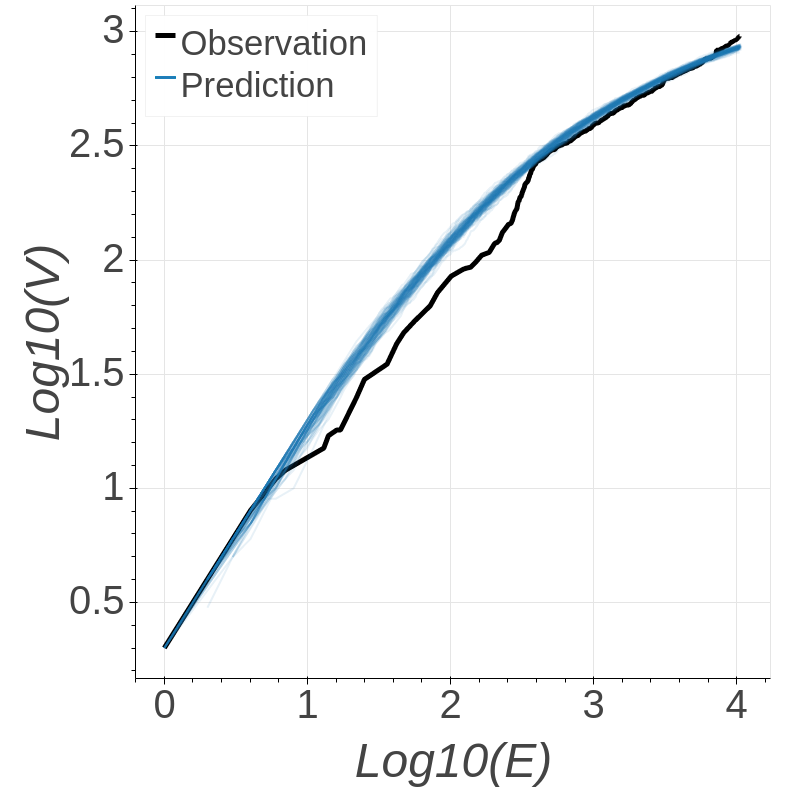}
		\caption{\# Vertices vs.~\# edges}\label{fig:fb_ev}
	\end{subfigure}
	\begin{subfigure}{0.32\textwidth}
		\centering\includegraphics[width=\textwidth]{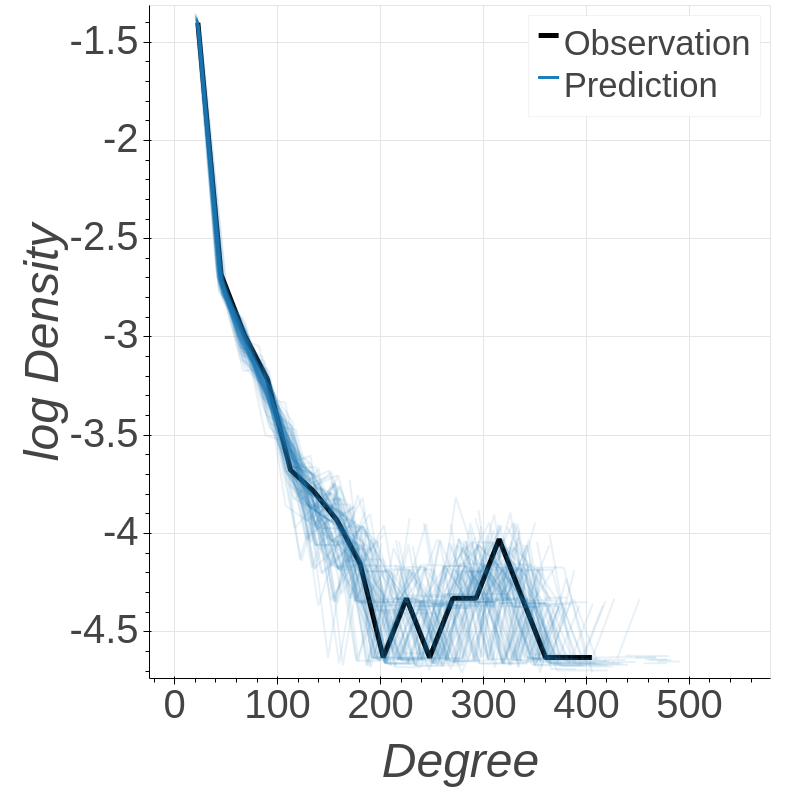}
		\caption{Degree density}\label{fig:fb_dd}
	\end{subfigure}
	\begin{subfigure}{0.32\textwidth}
		\centering\includegraphics[width=\textwidth]{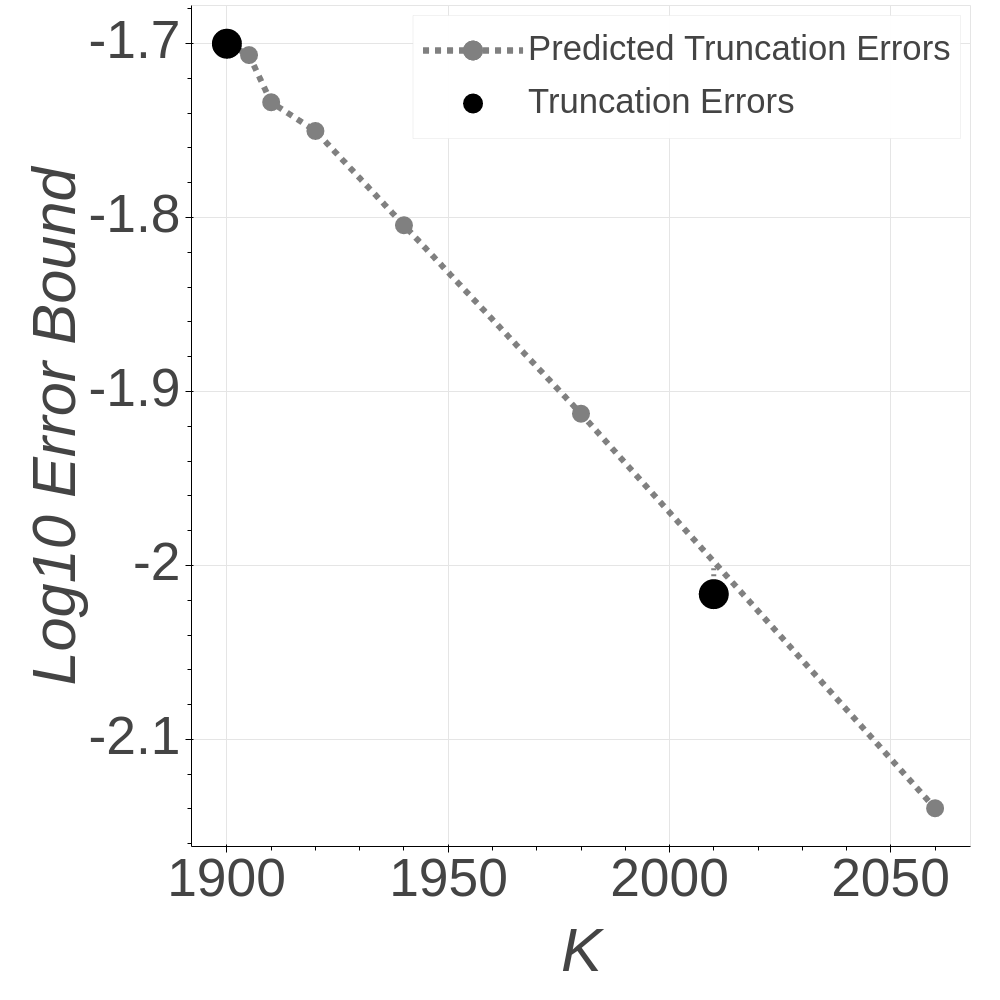}
		\caption{Truncation expansion}\label{fig:fb_exp}
	\end{subfigure}
	\caption{(\ref{fig:fb_ev},\ref{fig:fb_dd}): The characteristics of the observed Facebook-like network (black) and 100 samples from the posterior predictive distribution (blue). (\ref{fig:fb_exp}): The truncation adaptation process, with truncation level and error bound in each iteration (black circles) and predictive truncation errors (grey dashed lines).}\label{fig:fb_pre}
\end{figure}

\cref{fig:fb} shows the posterior marginal histograms for the hyperparameters $\alpha, \lambda, \gamma$
in both the first iteration (orange) and the second and final iteration (blue) of truncation adaptation.
The posterior distribution suggests that the network is dense (i.e.~$\alpha \approx 0$). This conclusion is 
supported both by the close match of 100 samples from the posterior predictive distribution, 
shown in \cref{fig:fb_ev,fig:fb_dd}, and the findings of past work using this data \cite{caron2017sparse}.
Further, as in the synthetic examples, the truncation adaptation terminates after two iterations; but
in this case the histograms do not change very much between the two. This is essentially because
the truncation error in the first iteration is relatively low ($\approx 0.02$), leading to 
a reasonably accurate truncated posterior and hence accurate predictions of the truncation
error at higher truncation levels.

\section{Conclusion}
In this paper, we developed methods for tractable generative simulation and
posterior inference in statistical models with discrete nonparametric priors
via finite truncation. We demonstrated that these approximate truncation-based
approaches are sound via theoretical error analysis.  
In the process, we also
showed that the nonzero rates of the (tractable) rejection representation of a
Poisson process are equal in distribution to the rates of the (intractable)
inverse L\'evy representation. Simulated and real network examples demonstrated
that the proposed methods are useful in selecting truncation levels 
for both forward generation and inference in practice.

\bibliographystyle{unsrt}
\bibliography{trunCRM}

\clearpage
\appendix
\section{Equivalence between nonzero rates from a rejection representation and the inverse L\'evy representation}
\bprfof{\cref{thm:levy}}
Denote $T_{k_1}$ be the first nonzero element that is generated
from the rejection representation from \cref{eq:rejecrep} and correspondingly,
denote $\Gamma_{k_1}$ be the jump of the unit-rate homogeneous Poisson process
on $\reals_+$ such that $T_{k_1} = \muinv(\Gamma_{k_1})$, where $\mu$ is the
proposal measure in the rejection representation. Let $f$ be a bounded 
continuous function. Then 
\[
\EE[f(T_{k_1})] &= \EE[f(\muinv(\Gamma_{k_1}))] \\
&= \EE\left[\sum_{j= 1}^{\infty}f(\muinv(\Gamma_j))\ind\left[\frac{\dee\nu}{\dee\mu}(\muinv(\Gamma_j))\ge U_j)\right]\prod_{i = 1}^{j-1}\ind\left[\frac{\dee\nu}{\dee\mu}(\muinv(\Gamma_i)) < U_i\right] \right]\\
&= \EE\left[\sum_{j=1}^{\infty}f(\muinv(\Gamma_j))\frac{\dee\nu}{\dee\mu}(\muinv(\Gamma_j))\prod_{i=1}^{j-1}\left(1-\frac{\dee\nu}{\dee\mu}(\muinv(\Gamma_i))\right) \right]\\
&= \sum_{j=1}^{\infty}\EE\left[f(\muinv(\Gamma_j))\frac{\dee\nu}{\dee\mu}(\muinv(\Gamma_j))\EE\left[\left.\prod_{i=1}^{j-1}\left(1-\frac{\dee\nu}{\dee\mu}(\muinv(\Gamma_i))\right)\right|\Gamma_j\right] \right]
\]
Note that given $\Gamma_j$, $\Gamma_i\distiid \distUnif(0, \Gamma_j)$, for $i = 1,\cdots, j-1$, so
\[
\EE\left[\left.\prod_{i=1}^{j-1}\left(1-\frac{\dee\nu}{\dee\mu}(\muinv(\Gamma_i))\right)\right|\Gamma_j\right] = \EE\left[1-\frac{\dee\nu}{\dee\mu}(\muinv(U)) \given \Gamma_j\right]^{j-1},
\]
where $U\distiid\distUnif(0, \Gamma_j)$.  Using the change of variable $y = \muinv(u)$, we obtain
\[
\EE\left[1-\frac{\dee\nu}{\dee\mu}(\muinv(U)) \given \Gamma_j\right] = 1 - \frac{1}{\Gamma_j}\int_{0}^{\Gamma_j}\frac{\dee\nu}{\dee\mu}(\muinv(u))\dee u 
= 1 - \frac{1}{\Gamma_j}\int_{\muinv(\Gamma_j)}^{\infty}\dee \nu.
\]
Therefore, using the same change of variable trick,
\[
\EE[f(T_{k_1})]  &=\sum_{j=1}^{\infty}\EE\left[f(\muinv(\Gamma_j))\frac{\dee\nu}{\dee\mu}(\muinv(\Gamma_j))\left(1-\frac{1}{\Gamma_j}\int_{\muinv(\Gamma_j)}^{\infty}\dee\nu \right)^{j-1} \right]\\
&= \sum_{j=1}^{\infty}\int_{0}^{\infty}f(\muinv(\gamma))\frac{\dee\nu}{\dee\mu}(\muinv(\gamma))\left(1-\frac{1}{\gamma}\int_{\muinv(\gamma)}^{\infty}\dee\nu\right)^{j-1}\frac{\gamma^{j-1}}{(j-1)!}e^{-\gamma}\dee \gamma\\
&=\sum_{j=1}^{\infty}\int_{0}^{\infty}f(y)\left(1-\mu[y, \infty)^{-1}\nu[y, \infty)\right)^{j-1}\frac{\mu[y, \infty)^{j-1}}{(j-1)!}e^{-\mu[y, \infty)}\nu(\dee y) \\
&=\int_{0}^{\infty}f(y)\sum_{j=1}^{\infty}\frac{\left(\mu[y, \infty)-\nu[y, \infty)\right)^{j-1}}{(j-1)!}e^{-\mu[y, \infty)}\nu(\dee y) \\
&= \int_{0}^{\infty}f(y)e^{-\nu[y, \infty)}\nu(\dee y).
\]
Suppose that $\theta_1$ is the first rate generated using the inverse L\'evy representation. Then
\[
\EE[f(\theta_1)] = \int_{0}^{\infty}f(\nu^{\leftarrow}(\gamma))e^{-\gamma}\dee\gamma.
\]
Making the change of variable $y = \nu^{\leftarrow}(\gamma)$, we obtain
\[
\EE[f(\theta_1)] = \int_{0}^{\infty}f(y)e^{-\nu[\gamma, \infty)}\nu(\dee y) = \EE[f(T_{k_1})].
\]
Therefore, the first nonzero rate $\theta_{k_1}$ from the rejection
representation has the same marginal distribution as the first rate
$\theta_1$ from the inverse L\'evy representation. 

We now employ an inductive argument. Suppose that we have shown
that the marginal distribution of first nonzero $M$ elements
$\Xi_{M}\defined (T_{k_1}, T_{k_2}, \cdots, T_{k_M})$ from
the rejection representation has the same marginal distribution as the first
$M$ elements $\Theta_M\defined (\theta_1, \cdots, \theta_M)$ from the inverse
L\'evy representation. To prove the same for $M+1$ elements, it 
suffices to show that the conditional distribution of
$T_{k_{M+1}}$ given $\Xi_M$ is equal to the conditional
distribution of $\theta_{M+1}$ given $\Theta_M$ when $\Xi_M = \Theta_M$. 

 Denote $\Gamma'_{j} = \sum_{i=1}^j e'_i$, where $e'_i\distiid \distExp(1)$,
and $U'_i \distiid \distUnif[0, 1]$. Then
\[
  &\EE[f(T_{k_{M+1}})|\Xi_M] = \EE[f(\muinv(\Gamma_{k_{M+1}}))|\Xi_M] = \EE\left[\left.\sum_{j=1}^{\infty}f(\muinv(\Gamma_{k_M}+\Gamma'_j))\ind[\dots]\right|\Xi_M \right],
\]
where $\ind[\dots]$ is shorthand for
\[
\ind[\dots] = \ind\left[\frac{\dee\nu}{\dee\mu}(\muinv(\Gamma_{k_M}+\Gamma'_j))\ge U'_j\right]\prod_{i=1}^{j-1}\ind\left[\frac{\dee\nu}{\dee\mu}(\muinv(\Gamma_{k_M}+\Gamma_i))<U'_i\right].
\]
Using steps similar to the base case,
\[
\EE[f(T_{k_{M+1}})|\Xi_M]
=\sum_{j=1}^{\infty}\EE\left[\left.f(\muinv(\Gamma_{k_M}+\Gamma'_j))\frac{\dee\nu}{\dee\mu}(\muinv(\Gamma_{k_M}+\Gamma'_j))\EE[\dots]^{j-1}\right|\Xi_M\right],
\]
where
\[
\EE[\dots] = \EE\left[\left.1-\frac{\dee\nu}{\dee\mu}(\muinv(\Gamma_{k_M}+U)) \right| \Gamma'_j, \Xi_M \right] \qquad U\dist \distUnif[0, \Gamma'_j].
\]
Making the change of variable $y = \muinv(\Gamma_{k_M} + u)$ as before, we obtain 
\[
\EE\left[\left.1-\frac{\dee\nu}{\dee\mu}(\muinv(\Gamma_{k_M}+U)) \right| \Gamma'_j,\Xi_M \right] &= 1 - \int_{0}^{\Gamma'_j}\frac{1}{\Gamma'_j}\frac{\dee\nu}{\dee\mu}(\muinv(\Gamma_{k_M}+u))\dee u \\
&= 1 - \frac{1}{\Gamma'_j}\int_{\muinv(\Gamma_{k_M}+\Gamma'_j)}^{\muinv(\Gamma_{k_M})}\dee\nu.
\]
Making another change of variables $y = \muinv(\Gamma_{k_M}+\gamma)$ in the original integral---and 
hence $\gamma = \mu[y,\infty) - \Gamma_{k_M} = \mu[y, \muinv(\Gamma_{k_M}))$---yields
\[
&\EE[f(T_{k_{M+1}})|\Xi_M]\\ =&\sum_{j=1}^{\infty}\int_{0}^{\infty}f(\muinv(\Gamma_{k_M}+\gamma))\frac{\dee\nu}{\dee\mu}(\muinv(\Gamma_{k_M}+\gamma))\EE[\dots]^{j-1}\frac{\gamma^{j-1}}{(j-1)!}e^{-\gamma}\dee\gamma \\
=&\sum_{j=1}^{\infty}\int_{0}^{\muinv(\Gamma_{k_M})}f(y)\left(1-\frac{1}{\gamma}\int_{y}^{\muinv(\Gamma_{k_M})}\dee\nu\right)^{j-1}\frac{\gamma^{j-1}}{(j-1)!}e^{-\gamma}\nu(\dee y)\\
=&\sum_{j=1}^{\infty}\int_{0}^{\muinv(\Gamma_{k_M})}f(y)\frac{(\mu[y, \muinv(\Gamma_{k_M})) - \nu[y, \muinv(\Gamma_{k_M})))^{j-1}}{(j-1)!}e^{-\mu[y, \muinv(\Gamma_{k_M}))}\nu(\dee y)\\
=&\int_{0}^{\muinv(\Gamma_{k_M})}f(y)e^{-\nu[y,\muinv(\Gamma_{k_M}))}\nu(\dee y)
=\int_{0}^{T_{k_M}}f(y)e^{-\nu[y,T_{k_M})}\nu(\dee y).
\]
On the other hand,
\[
&\EE[f(\theta_{M+1})|\Theta_M] = \EE[f(\nu^{\leftarrow}(\Gamma_M + \Gamma'_1))|\Theta_{M}] 
= \int_{0}^{\infty}f(\nu^{\leftarrow}(\Gamma_M + \gamma))e^{-\gamma}\dee\gamma\\
=& \int_{0}^{\nu^{\leftarrow}(\Gamma_M)}f(y)e^{-\nu[y,\nu^{\leftarrow}(\Gamma_M))}\nu(\dee y)
= \int_{0}^{\theta_M}f(y)e^{-\nu[y, \theta_M)}\nu(\dee y).
\]
Thus the distribution of the $(M+1)^{\text{th}}$ nonzero rate in the rejection representation $T_{k_{M+1}}$ given $\Xi_M$ is equal to the  distribution of the $(M+1)^{\text{th}}$ rate from the inverse L\'evy representation $\theta_{M+1}$ given $\Theta_M$ when $\Xi_M = \Theta_M$. 
\eprfof

\section{Truncation error bounds for self-product measures}\label{sec:truncproofs}
\bprfof{\cref{lem:marginal}}
Denote $\tilde{\Theta} = \{\Theta^{(d)}, \Theta_K^{(d)} \}$. 
Denote the marginal probability mass function (PMF) of $E_N\in\mcN_d$ and
$E_{N,K}\in\mcN_d$ as $P_N$ and $P_{N,K}$, and denote their PMFs given
$\tilde{\Theta}$ as $f(x|\tilde{\Theta})$ and $f_{K}(x|\tilde{\Theta})$
respectively.  
\[
& \tvd{P_{N}}{P_{N,K}}\\
&= \frac{1}{2}\sum_{x\in\mcN_d}\Big| P_{N}(x)-P_{N,K}(x)\Big| \\
&= \frac{1}{2}\sum_{x\in\mcN_d}\left|\EE[f(x|\tilde{\Theta})] - \EE[f_K(x|\tilde{\Theta})]\right| \\
&\leq \frac{1}{2}\Pr(I_N\leq K)\sum_{x\in \mcN_d}\left|\EE[f(x|\tilde{\Theta})|I_N\leq K] - \EE[f_K(x|\tilde{\Theta})|I_N\leq K]\right| \\
&\,\,\, + \frac{1}{2}\Pr(I_N> K)\sum_{x\in\mcN_d}\left|\EE[f(x|\tilde{\Theta})|I_N> K] - \EE[f_K(x|\tilde{\Theta})|I_N> K]\right| \\
\]
Conditioned on $I_N\leq K$, $f(x|\tilde{\Theta}) = f_K(x|\tilde{\Theta})$ under both the independent and categorical likelihood. 
So
\[
\EE[f(x|\tilde{\Theta})|I_N\leq K] - \EE[f_K(x|\tilde{\Theta})|I_N\leq K] = 0.
\]
By Fubini’s Theorem,
\[
& \tvd{P_{N}}{P_{N,K}}\\
\leq& \frac{1}{2}\Pr(I_N>K)\sum_{x\in\mcN_d}\left(\ \EE[f(x|\tilde{\Theta})|I_N>K] + \EE[f_K(x|\tilde{\Theta})|I_N>K]\right) \\
=& \frac{1}{2}\Pr(I_N>K)\left(\EE\left[\left.\sum_{x\in\mcN_d} f(x|\tilde{\Theta}) \right|I_N>K\right] + \EE\left[\left. \sum_{x\in\mcN_d} f_K(x |\tilde{\Theta}) \right| I_N>K\right] \right)\\
=& \Pr(I_N>K) = 1-\Pr\left(I_N\le K\right).
\]
\eprfof

\bprfof{\cref{thm:truncerror,thm:conditionaltrunc}}
Denote the set of indices
\[
\mcI_{\ell, K} \defined
 \{i\in \nats_{\neq}^d: 1 \le i_1, \cdots, i_\ell\le K,\  K+1\le i_{\ell+1},\cdots, i_d < \infty \}
\]
such that $i\in \mcI_{\ell, K}$ indicates that the first $\ell$ elements of $i$ belong to 
the truncation, and the remaining $d-\ell$ elements belong to the tail. 
By Jensen's inequality,
\[
&\Pr\left(I_N \le K \given U_{1:K}, \Gamma_{1:K}\right)\nonumber\\
&=\EE\left[\exp\left\{N\sum_{\ell=0}^{d-1}\binom{d}{\ell}\sum_{i\in \mcI_{\ell, K}}\log\pi\left(\vartheta_i\right) \right\} \given U_{1:K}, \Gamma_{1:K}\right] \nonumber\\
&\geq \exp\left\{N\sum_{\ell=0}^{d-1}\binom{d}{\ell}\EE\left[\sum_{i\in\mcI_{\ell, K}}\log\pi\left(\vartheta_i\right)\given U_{1:K}, \Gamma_{1:K}\right] \right\}. \label{eq:ddimineq}
\]
This equation arises by noting that $I_N \leq K$ if and only if for all $i$
involving an index $i_j > K$, the count of edge $i$ is 0 after $N$ rounds; 
the factor $\binom{d}{\ell}$ accounts for the fact that
$\vartheta_i = \prod_{j=1}^d \theta_{i_j}$ is independent of the ordering
of the $i_j$.

Consider a single term $\sum_{i\in\mcI_{\ell,K}}\log\pi(\vartheta_i)$ in the above sum.
Since we are conditioning on $U_{1:K}, \Gamma_{1:K}$, we have that $\theta_{i_1},\dots,\theta_{i_\ell}$ 
are fixed in the expectation, and the remaining steps $\Gamma_{K+1},
\Gamma_{K+2}, \cdots$ are the ordered jumps of a unit rate homogeneous Poisson process on $\left[\Gamma_K, \infty\right)$. 
By the marking property of the Poisson process
\cite{kingman1992poisson}, conditioned on $\Gamma_K$, we have that $(U_i, \Gamma_i)_{i=K+1}^{\infty}$ 
is a Poisson process on $\reals_+\times\left[\Gamma_K, \infty\right)$ with
rate measure $g(\dee u)\dee \gamma$.
Thus we apply the Slivnyak-Mecke theorem \cite{last2017lectures} to the remaining $d-\ell$ indices to obtain
\[
 &\EE\left[\left.\sum_{i\in\mcI_{\ell, K}}\log\pi\left(\vartheta_i\right)\right|U_{1:K},\Gamma_{1:K}\right]\\
 &\EE\left[\left.\sum_{i\in\mcI_{\ell, K}}\log\pi\left(\prod_{j=1}^{\ell}\theta_{i_j}\prod_{j=\ell+1}^{d}\tau(U_{i_j}, \Gamma_{i_j})\right)\right|U_{1:K},\Gamma_{1:K}\right]\\
=&\sum_{i_1\ne\cdots\ne i_\ell \leq K}\EE\left[\int_{\reals_+^{d-\ell}}\log\pi\left(\prod_{j=1}^{\ell}\theta_{i_j}\prod_{j=\ell+1}^d\tau(U_{i_j}, \Gamma_K+\gamma_{j})\right)\prod_{j=\ell+1}^{d} \dee\gamma_j \given U_{1:K}, \Gamma_{1:K}\right]\\
=&\sum_{\begin{subarray}{c}\mcL \subseteq [K]\\|\mcL| = \ell\end{subarray}}\int_{[\Gamma_K, \infty)^{d-\ell}}\!\!\!\!\EE\left[\log\pi\left(\prod_{j\in\mcL}\theta_{j}\cdot\prod_{j=1}^{d-\ell}\tau(U'_{j}, \gamma_{j})\right) \given \theta_{1:K}\right]\dee\gamma_{1:d-\ell},
\]
where $U'_j \distiid g$. Substitution of this expression into \cref{eq:ddimineq} yields the 
result of \cref{thm:conditionaltrunc}.
Next, we consider the bound on the marginal probability $\Pr\left(I_N\leq K\right)$.
By Jensen's inequality applied to \cref{eq:ddimineq} and following the previous derivation, we 
find that
\[
\Pr\left(I_N \le K \right) \geq \exp\left\{N\sum_{\ell=0}^{d-1}\binom{d}{\ell}\EE\left[\sum_{i\in\mcI_{\ell, K}}\log\pi\left(\vartheta_i\right)\right] \right\},
\]
and
\[
&\EE\left[\sum_{i\in\mcI_{\ell, K}}\log\pi\left(\vartheta_i\right)\right]\\
=&\sum_{i_1\ne\cdots\ne i_\ell \leq K}\EE\left[\int_{[\Gamma_K, \infty)^{d-\ell}}\!\!\!\!\log\pi\left(\prod_{j=1}^{\ell}\theta_{i_j}\prod_{j=1}^{d-\ell}\tau(U'_{j}, \gamma_{j})\right)\dee\gamma_{1:d-\ell}\right].
\]
Using the fact that conditioned on $\Gamma_K$, $\Gamma_{1:K-1}$ are uniformly distributed on $[0, \Gamma_K]$,
and that at most one $i_j$ can be equal to $K$,
\[
&\sum_{i_1\ne\cdots\ne i_\ell \leq K}\EE\left[\int_{[\Gamma_K, \infty)^{d-\ell}}\!\!\!\!\log\pi\left(\prod_{j=1}^{\ell}\theta_{i_j}\prod_{j=1}^{d-\ell}\tau(U'_{j}, \gamma_{j})\right)\dee\gamma_{1:d-\ell}\right]\\
&=\frac{(K-1)!}{(K-1-\ell)!}\EE\left[\Gamma_K^{-\ell}\int_{[0,\Gamma_K]^\ell\times[\Gamma_K, \infty)^{d-\ell}}\!\!\!\!\log\pi\left(\prod_{j=1}^{d}\tau(U_{j}, \gamma_{j})\right)\dee\gamma_{1:d}\right]\\
&+\ell\frac{(K-1)!}{(K-\ell)!}\EE\left[\Gamma_K^{-(\ell-1)}\int_{[0,\Gamma_K]^\ell\times[\Gamma_K, \infty)^{d-\ell}}\!\!\!\!\delta_{\gamma_1=\Gamma_K}\log\pi\left(\prod_{j=1}^{d}\tau(U_{j}, \gamma_{j})\right)\dee\gamma_{1:d}\right],
\]
where the first and second terms arise from portions of the sum where all indices satisfy $i_j \neq K$ 
and one index satisfies $i_j = K$, respectively.

To complete the result, we study the asymptotics of the marginal probability bound.
It follows from \cref{eq:hcond} that 
$I_N < \infty$ almost surely. Therefore 
\[
\lim_{K\rightarrow\infty}\EE\left[\exp\left\{N\sum_{\ell=0}^{d-1}\binom{d}{\ell}\sum_{i\in\mcI_{\ell, K}}\log\pi\left(\vartheta_i\right) \right\}\right] 
=\lim_{K\rightarrow\infty} \Pr(I_N\le K) = 1.
\]
It then follows from \cite[ Lemma B.1]{campbell2019truncated} and continuous mapping theorem that 
\[
\sum_{\ell=0}^{d-1}\left[\binom{d}{\ell}\sum_{i\in\mcI_{\ell, K}}\log\pi\left(\vartheta_i\right)\right] \convp 0\qquad\text{as}\,\,K\to\infty.
\]
Since this sequence is monotonically increasing in $K$, we have that
\[
B_K = \EE\left[\sum_{\ell=0}^{d-1}\left[\binom{d}{\ell}\sum_{i\in\mcI_{\ell, K}}-\log\pi\left(\vartheta_i\right)\right]\right] \to 0\qquad\text{as}\,\,K\to\infty.
\]

\eprfof

\subsection{Proof of \cref{thm:truncnorm}}\label{sec:prftruncnorm}
We first state an useful results which states that if one perturbs the probabilities of a countable discrete distribution by i.i.d. $\distGumbel(0, 1)$ random variables, the arg max of the resulting set is a sample from that distribution.
\bnlem{\cite[Lemma 5.2]{campbell2019truncated}}
Let $(p_j)_{j=1}^{\infty}$ be a countable collection of positive numbers such that $\sum_j p_j<\infty$ and let $\bar{p}_j = \frac{p_j}{\sum_k p_k}$. If $(W_j)_{j=1}^{\infty}$ are i.i.d $\distGumbel(0, 1)$ random variables, then $\argmax_{j\in \nats} W_j + \log p_j$ exists, is unique a.s., and has distribution
\[
\argmax_{j\in \nats} W_j + \log p_j \sim \distCat\left((\bar{p}_j)_{j=1}^{\infty}\right).
\]
\enlem
\bprfof{\cref{thm:truncnorm}}
Since the $N$ edges from the categorical likelihood process are \iid categorical random variables, by Jensen's inequality,
\[
\Pr(I_N\le K) = \EE\left[\Pr(I_N\le K|\Theta)\right] = \EE\left[\Pr(I_1\le K|\Theta)^N\right] \ge \EE\left[\Pr(I_1\le K|\Theta)\right]^N.
\]
Next, since $\vartheta_{i} = \prod_{j=1}^{d}\theta_{i_j}$, we can 
simulate a categorical
variable with probabilities proportional to $\vartheta_i$, $i\in
\nats_{\neq}^d$ by sampling $d$ indices $(J_1, \cdots, J_d)$
independently from a categorical distribution with probabilities proportional to
$(\theta_1, \theta_2, \dots)$ and discarding samples where
$J_j = J_k$ for some $1\le j, k\le d$. Denote $\theta'_k = \theta_k/\sum_k\theta_k$ to be the
normalized rates, $P_{J, K}\defined \sum_{j=J}^{K}
\theta'_j$, and the event $\mcQ \defined \{J_j\ne J_k, \forall 1\le j, k\le d
\}$. Then 
\[
\Pr(I_1\le K|\Theta) = \Pr(1\le J_j \le K,\ 1\le j\le d\ |\ \mcQ, \Theta).
\]
Since  the normalized rates $\theta'_k$ are generated from the inverse L\'evy representation,
they are monotone decreasing. Therefore
\[
\Pr(1\le J_j \le K,\ 1\le j\le d\ |\ \mcQ, \Theta) &\ge \frac{P_{1,K}}{1 - 0}\cdot\frac{P_{2, K}}{1-P_{1, 1}}\cdots\frac{P_{d, K}}{1-P_{1, d-1}}\\ 
&\ge \left(\frac{P_{d, K}}{1-P_{1, d-1}}\right)^d.
\]
By Jensen's inequality,
\[
\EE[\Pr(1\le J_j \le K,\ 1\le j\le d\ |\ \mcQ, \Theta)]\ge \EE\left[\frac{P_{d, K}}{1-P_{1, d-1}}\right]^d.
\]
Note that for a categorical random variable $J$ with class probabilities given
by $\theta'_j/(1-P_{1, d-1})$, $j\ge d$, the quantity $P_{d,
K}/(1-P_{1, d-1})$ is the probability that $d\le J\le K$.  So by the infinite
Gumbel-max sampling lemma, 
\[
\EE\left[\frac{P_{d,K}}{1-P_{1, d-1}}\right] = \Pr\left(d\le \argmax_{j\ge d} (\log \theta_j + W_j) \le K\right), \quad W_j\distiid \distGumbel(0, 1),
\]
where we can replace $\theta'_j$ with the unnormalized $\theta_j$ because the normalization does not affect the $\argmax$.
Denoting
\[
M_K\defined \max_{d \le k \le K} \log \nuinv(\Gamma_k) + W_k, & &
M_{K+}\defined \sup_{k>K}\log\nuinv(\Gamma_k) + W_k,
\]
we have that
\[
\Pr(I_N \leq K) \ge \left(1 - \EE[\Pr(M_K < M_{K+} | \Gamma_K)]\right)^{N\cdot d},
\]
and so the remainder of the proof focuses on the conditional expectation.
Conditioned on $\Gamma_K$,
\[
M_K \eqd \max\left\{\log\nuinv(\Gamma_K) + W_K, \ \max_{d \le k \le K}\log\nuinv(u_k) + W_k \right\}.
\]
The cumulative distribution function and the probability density function of the Gumbel distribution $\distGumbel(0, 1)$ is 
\[
F(x) = e^{-e^{-x}}, \quad f(x) = e^{-(x + e^{-x})}.
\]
So
\[
\Pr(\log\nuinv(\Gamma_K) + W_K\le x\ |\ \Gamma_K) = e^{-\nuinv(\Gamma_K)e^{-x}},
\]
and
\[
\Pr(\log\nuinv(u_k) + W_k \le x\ |\ \Gamma_K) = \int_{0}^{1}e^{-\nuinv(\Gamma_K u)e^{-x}}\dee u.
\]
Therefore,
\[
\Pr(M_K\le x\ |\ \Gamma_K) = \left(\int_{0}^{1}e^{-\nuinv(\Gamma_K u)e^{-x}}\dee u\right)^{K-d}e^{-\nuinv(\Gamma_K)e^{-x}}.
\]
Denote $Q(u, t) = e^{-\nuinv(u)e^{-t}}$, and 
\[
\Pr(M_K\le x\ |\ \Gamma_K) = \left(\int_{0}^{1}Q(\Gamma_K u, x)\dee u\right)^{K-d}Q(\Gamma_K, x).
\]
Conditioned on $\Gamma_K$, the tail $\Gamma_{K+1}, \Gamma_{K+2}, \cdots$ is a unit rate homogeneous Poisson process on $[\Gamma_K, \infty)$ that is independent of $\Gamma_{1},\cdots, \Gamma_{K-1}$. So conditioned on $\Gamma_K$,
\[
M_{K+} \eqd \sup_{k\ge 1}\log\nuinv(\Gamma_K + \Gamma_k') + W_k,
\]
where $\Gamma_k'$ is unit rate of homogeneous Poisson process on $\reals_+$. Since $\Gamma_k'$ is a Poisson process, $\log\nuinv(\Gamma_K + \Gamma_k') + W_k$ is also a Poisson process with the rate measure 
\[
\left(\int_{0}^{\infty}e^{-(t-\log\nuinv(\Gamma_K + \gamma))-e^{-(t-\log\nuinv(\Gamma_K + \gamma))}}\dee\gamma \right)\dee t.
\]
$\Pr(M_{K+}\le x\ |\ \Gamma_K)$ is the probability that no atom of the Poisson process is greater than $x$. For a Poisson process with rate measure $\mu(\dee t)$, this probability is $e^{-\int_{x}^{\infty}\mu(\dee t)}$,
\[
\Pr(M_{K+}\le x \ | \ \Gamma_K) &= e^{-\int_{x}^{\infty}\left(\int_{0}^{\infty}e^{-(t-\log\nuinv(\Gamma_K+\gamma)) - e^{-(t - \log\nuinv(\Gamma+\gamma))}}\dee\gamma \right) \dee t}\\
&=e^{\int_{0}^{\infty}\left(e^{-\nuinv(\Gamma_K + \gamma)e^{-x}} - 1 \right)\dee \gamma}\\
&=e^{\int_{0}^{\infty}Q(\Gamma_K+\gamma,\, x) - 1\, \dee\gamma},
\]
where the second equation comes from the fact that the inner integrand is the probability density function of a Gumbel distribution.
Therefore,
\[
 &\Pr(M_K< M_{K+}\ |\ \Gamma_K) \\
=& \int_{-\infty}^{\infty}\Pr(M_K\le x\ | \ \Gamma_K)\frac{\dee}{\dee x}\Pr(M_{K+}\le x\ |\ \Gamma_K)\dee x \\
=& \int_{-\infty}^{\infty}Q(\Gamma_K, x)\left(\int_{0}^{1}Q(\Gamma_K u, x)\dee u\right)^{K-d}\left(\frac{\dee}{\dee x}e^{\int_{0}^{\infty}Q(\Gamma_K+\gamma,\, x) - 1\, \dee\gamma}\right) \dee x.
\]
For the categorical variable $J$ with class probabilities given by $\theta'_j/(1-P_{1, d-1})$, $j\ge d$, it holds that $\Pr(d\le J\le K) \uparrow 1$ as $K\to \infty$. By the monotone convergence theorem
\[
B_K = \Pr\left(\argmax_{j\ge d} (\log \theta_j + W_j) > K \right)  = 1-\EE\left[\frac{P_{d,K}}{1-P_{1,d}} \right]\to 0.
\]

\eprfof

\section{Error bounds for edge-exchangeable networks}
\subsection{Rejection representation}
We first derive the specific form of $B_{K}$ in \cref{cor:1} for the rejection representation from \cref{eq:rejecrep}, where
\[
\tau(U, \Gamma) = \muinv(\Gamma)\mathds{1}\left[\frac{\dee \nu}{\dee\mu}(\muinv(\Gamma))\geq U\right],
\]
and $U\sim\distUnif[0, 1]$. So in \cref{cor:1},
\[
& B_{K,1} \\
= &\EE\left[\int_{\reals_+^4}-\log\pi\left(\tau(u_1, \gamma_1+\Gamma_K)\tau(u_2, \gamma_2+\Gamma_K)\right)\mathrm{d}\gamma_1\mathrm{d}\gamma_2 g(u_1)\mathrm{d}u_1g(u_2)\mathrm{d}u_2\right] \\
=& \EE\left[\int_{\Gamma_K}^{\infty}\int_{\Gamma_K}^{\infty}\int_{0}^{1}\int_{0}^{1}-\log\pi\left(\tau(u_1, \gamma_1)\tau(u_2, \gamma_2)\right)\mathrm{d}u_1\mathrm{d}u_2\mathrm{d}\gamma_1\mathrm{d}\gamma_2\right].
\]
Since $\pi(0) = 1$, $\log\pi(0) = 0$, we can take the indicator in $\tau$ out 
of the function  $\log\pi(\tau(u_1, \gamma_1)\tau(u_2, \gamma_2))$ to obtain
\[
&\int_{0}^{1}\int_{0}^{1}\log\pi\left(\tau(u_1, \gamma_1)\tau(u_2, \gamma_2)\right)\mathrm{d}u_1\mathrm{d}u_2\\
=&\int_{0}^{1}\int_{0}^{1}\log\pi\left(\muinv(\gamma_1)\muinv(\gamma_2)\right)\mathds{1}\left[\frac{\dee\nu}{\dee\mu}(\muinv(\gamma_1))\geq u_1\right]\mathds{1}\left[\frac{\dee\nu}{\dee\mu}(\muinv(\gamma_2))\geq u_2\right]\mathrm{d}u_1\mathrm{d}u_2\\
=& \log\pi\left(\muinv(\gamma_1)\muinv(\gamma_2)\right)\frac{\dee\nu}{\dee\mu}(\muinv(\gamma_1))\frac{\dee\nu}{\dee\mu}(\muinv(\gamma_2)).
\]
Transforming variables via 
$x_1 = \muinv(\gamma_1)$ and $x_2 = \muinv(\gamma_2)$, and noting 
that $\muinv(\Gamma_K)\ge x \iff \Gamma_K\le \mu[x, \infty)$, 
\[
B_{K,1}
&= \EE\left[\int_{0}^{\muinv(\Gamma_K)}\int_{0}^{\muinv(\Gamma_K)}-\log\pi(x_1x_2)\frac{\mathrm{d}\nu}{\mathrm{d}\mu}(x_1)\frac{\dee\nu}{\dee\mu}(x_2)\mu(\mathrm{d}x_1)\mu(\mathrm{d}x_2)\right] \nonumber\\
&=\EE\left[\int_{\reals_+^2}-\log\pi(x_1x_2)\mathds{1}[x_1\leq \muinv(\Gamma_K)]\mathds{1}[x_2\leq\muinv(\Gamma_K)]\nu(\mathrm{d}x_1)\nu(\mathrm{d}x_2)\right] \nonumber\\
&= \int_{\reals_+^2}-\log\pi(x_1x_2)\EE\left[\mathds{1}[\Gamma_K\leq\mu[\max\{x_1, x_2\},  \infty)]\right]\nu(\dee x_1)\nu(\dee x_2) \nonumber\\
&= \int_{\reals_+^2}-\log\pi(x_1x_2)F_K\left(\mu[\max\{x_1,x_2\}, \infty)\right)\nu(\dee x_1)\nu(\dee x_2), \label{eq:bnk1}
\]
where $F_K$ is the cumulative distribution function of $\Gamma_K$. Using the same variable transformation again,
\[
&\frac{1}{2(K-1)}B_{K,2} \\
=& \EE\left[\int_{0}^{\Gamma_K}\int_{\reals_+^3}-\frac{1}{\Gamma_K}\log\pi\left(\tau(u_1, \gamma_1)\tau(u_2, \gamma_2+\Gamma_K)\right)g(u_1)\mathrm{d}u_1g(u_2)\mathrm{d}u_2 \mathrm{d}\gamma_1\mathrm{d}\gamma_2\right] \\
=& \EE\left[\int_{\muinv(\Gamma_K)}^{\infty}\int_{0}^{\muinv(\Gamma_K)}-\frac{1}{\Gamma_K}\log\pi(x_1x_2)\frac{\dee\nu}{\dee\mu}(x_1)\frac{\mathrm{d}\nu}{\mathrm{d}\mu}(x_2)\mu(\mathrm{d}x_1)\mu(\mathrm{d}x_2)\right] \\
=& \int_{\reals_+^2}-\EE\left[\frac{1}{\Gamma_K}\mathds{1}[\mu[x_1,\infty)\leq \Gamma_K\leq\mu[x_2,\infty)]\right]\log\pi(x_1x_2)\nu(\mathrm{d}x_1)\nu(\mathrm{d}x_2).
\]
\[
B_{K,3} & = 2\EE\left[\int_{\reals_+^3}-\log\pi\left(\tau(u_1, \Gamma_K)\tau(u_2,\gamma_2 + \Gamma_K)\right)g(u_1)\dee u_1 g(u_2) \dee u_2\dee \gamma_2\right] \\
& = 2\EE\left[\int_{\Gamma_K}^{\infty}-\log\pi\left(\muinv(\Gamma_K)\muinv(\gamma_2)\right)\frac{\dee\nu}{\dee\mu}(\muinv(\Gamma_K))\frac{\dee\nu}{\dee\mu}(\muinv(\gamma_2))\dee \gamma_2\right] \\
& = 2\EE\left[\int_{0}^{\muinv(\Gamma_K)}-\log\pi\left(\muinv(\Gamma_K)x\right)\frac{\dee\nu}{\dee\mu}(\muinv(\Gamma_K))\frac{\dee\nu}{\dee\mu}(x)\mu(\dee x)\right] \\
& = 2\EE\left[\int_{\reals_+}-\log\pi\left(\muinv(\Gamma_K)x\right)\frac{\dee\nu}{\dee\mu}(\muinv(\Gamma_K))\mathds{1}[\Gamma_K\leq\mu[x,\infty)] \nu(\dee x)\right].
\]

\subsection{Beta-independent Bernoulli process network}\label{pf:beta}
\paragraph{Dense network} When $\alpha=0$, 
\[
\nu(\dee\theta) = \gamma\lambda\theta^{-1}(1-\theta)^{\lambda-1}\dee\theta, \quad \mu(\dee\theta) = \gamma\lambda\theta^{-1}\dee\theta,
\] 
and
\[ 
\frac{\dee\nu}{\dee\mu}=(1-\theta)^{\lambda-1},\quad \mu\left([x, 1]\right) = -\gamma'\log x,\quad\muinv(u) = e^{-u/\gamma'}. 
\]
Substituting $\nu(\mathrm{d}x), \mu(\mathrm{d}x)$ and $\pi(\theta)$ into \cref{eq:bnk1} and 
noting that the integrand is symmetric around the line $x_1=x_2$, 
\[
B_{K,1} &= 2\int_{0}^{1}\int_{x_1}^{1}-\log(1-x_1x_2)F_K\left(\mu[x_2, 1]\right)\gamma\lambda x_2^{-1}(1-x_2)^{\lambda-1}\dee x_2\nu(\dee x_1).
\]
Next, note that $0\le F_K\left(\mu[x_2, 1] \right)\le F_K\left(\mu[x_1, 1] \right)$ when $x_2 \ge x_1$ and 
\[
0\ge\log(1-x_1x_2)\ge -\left(\frac{1}{1-x_1x_2}-1\right) = -\left(\frac{x_1x_2}{1-x_1x_2}\right)\ge -\left(\frac{x_1x_2}{1-x_1}\right).
\]
So
\[
B_{K,1} &\le 2\gamma\lambda\int_{0}^{1}F_K\left(\mu([x_1, 1])\right)\int_{x_1}^{1}\frac{x_1x_2}{1-x_1}x_2^{-1}(1-x_2)^{\lambda-1}\dee x_2\nu(\dee x_1) \\
&= 2\gamma^2\lambda\int_{0}^{1}(1-x)^{2(\lambda-1)}F_K\left(\mu([x, 1])\right)\dee x.
\]
For any $a\in (0, 1)$, dividing the integral into two parts and bounding each part separately,  
\[
B_{K,1}\le 2\gamma^2\lambda\left[\int_{0}^{a}(1-x)^{2(\lambda-1)}\dee x + \int_{a}^{1}F_K\left(-\gamma'\log x\right)(1-x)^{2(\lambda-1)}\dee x \right].
\]
Assume for the moment that $\lambda\ne 0.5$. Use the fact that $F_K(t)\le (3t/K)^K$ and note also that when $x\in [a, 1]$, $-\log x \le \left[-\log a/(1-a)\right](1-x)$,
\[
& \frac{1}{2\gamma^2\lambda}B_{K,1}\\
\le& \frac{1-(1-a)^{2\lambda-1}}{2\lambda-1} + \int_{a}^{1}\left(\frac{3\gamma'\left[-\log a/(1-a)\right](1-x)}{K} \right)^K(1-x)^{2\lambda-2}\dee x  \\
=& \frac{1-(1-a)^{2\lambda-1}}{2\lambda-1} + \left(-\frac{3\gamma'\log a}{K}\right)^K\frac{1}{K+2\lambda-1}(1-a)^{2\lambda-1}.
\]
It can be seen that there exists a constant value $c>1$ such that $3\gamma'\log c = 1/c$. Setting $a = c^{-K}$ and using the first order Taylor's expansion to approximate the first term,  it can be seen that
\[
\frac{1-(1-a)^{2\lambda-1}}{2\lambda-1} \sim c^{-K}, \quad \left(-\frac{3\gamma'\log a}{K}\right)^K\frac{1}{K+2\lambda-1}(1-a)^{2\lambda-1}\sim \frac{1}{K}c^{-K}.
\]
Therefore, there exists $K_0$ such that when $K > K_0$,  
\[
B_{K,1} \le 4\gamma^2\lambda c^{-K}. 
\] 
If $\lambda=0.5$,
\[
B_{K,1} &\le 2\gamma^2\lambda\int_{0}^{1}(1-x)^{-1}F_K\left(\mu([x,1])\right)\dee x \\
&\le 2\gamma^2\lambda\left[-\log(1-a) + \left(\frac{3\gamma'}{K}\right)^K\frac{1}{K}(-\log a)^K\right],
\]
which can be bounded similarly by choosing the same $c$ and setting $a = c^{-K}$. Therefore we can find a constant $K_0$ and for $K > K_0$,
\[
B_{K,1} \le 4\gamma^2\lambda c^{-K}.
\] 
Next, the term $B_{K,2}$ is bounded via
\[
&\frac{1}{2}B_{K,2}\\
=& \int_{0}^{1}\int_{x_1}^{1}-\log(1-x_1x_2)\left[F_{K-1}\left(\mu([x_1, 1])\right)  - F_{K-1}\left(\mu([x_2, 1])\right)\right]\nu(\dee x_2)\nu(\dee x_1) \\
\le&\int_{0}^{1}\int_{x_1}^{1}-\log(1-x_1x_2)F_{K-1}\left(\mu[x_1, 1]\right)\nu(\dee x_2)\nu(\dee x_1).
\]
This has exactly the same form as $B_{K,1}$, except that the CDF of $\Gamma_K$ is 
replaced with that of $\Gamma_{K-1}$. Therefore, it can be shown that for large $K$,
\[
B_{K,2} \le 4\gamma^2\lambda c^{-(K-1)}. 
\] 
Finally, $B_{K,3}$ may be expressed as 
\[
&\frac{1}{2}B_{K,3}\\
=& \EE\left[\int_{0}^{1}-\log\left(1-\muinv(\Gamma_K)x\right)\frac{\dee \nu}{\dee\mu}\left(\muinv(\Gamma_K)\right)\mathds{1}\left[x\le\muinv(\Gamma_K)\right]\right]\nu(\dee x) \\
=& \gamma\lambda\EE\left[\int_{0}^{\muinv(\Gamma_k)}-\log\left(1-\muinv(\Gamma_K)x\right)\left(1-\muinv(\Gamma_K)\right)^{\lambda-1}x^{-1}(1-x)^{\lambda-1}\dee x \right].
\]
Since $\log\left(1-\muinv(\Gamma_K)x\right) \ge  - \muinv(\Gamma_K)x / \left(1-\muinv(\Gamma_K)\right)$,
\[
B_{K,3}&\le 2\gamma\lambda\EE\left[\left(1-\muinv(\Gamma_K)\right)^{\lambda-2}\muinv(\Gamma_K)\int_{0}^{\muinv(\Gamma_K)}(1-x)^{\lambda-1}\dee x \right] \\
&\ge 2\gamma\EE\left[\left(1-\muinv(\Gamma_K)\right)^{\lambda-2}\muinv(\Gamma_K) \right]\\
&= 2\gamma\int_{0}^{\infty}\left(1-e^{-x/\gamma'}\right)^{\lambda-2}e^{-x/\gamma'}\ \frac{x^{K-1}}{\Gamma(K)}e^{-x}\dee x.
\]
We split the analysis of this term into two cases. In the first case, assuming $\lambda\ge2$, we have that 
\[
B_{K,3} \le 2\gamma\int_{0}^{\infty}\frac{x^{K-1}}{\Gamma(K)}e^{-x(1+\gamma')/\gamma'}\dee x = 2\gamma\left(\frac{\gamma'}{1+\gamma'}\right)^K.
\]
On the other hand, if $\lambda<2$, we bound the integral over $[0, \gamma']$ and over $[\gamma', \infty)$ separately. Since $1-e^{-x}\ge x^2$ for $x\in [0, 1]$,
\[
& \frac{1}{2\gamma}B_{K,3}\\
\le& \int_{0}^{\gamma'}\left(\frac{x}{\gamma'}\right)^{2(\lambda-2)}\frac{x^{K-1}}{\Gamma(K)}e^{-x(1+\gamma')/\gamma'}\dee x + (1-e^{-1})^{\lambda-2}\int_{\gamma'}^{\infty}\frac{x^{K-1}}{\Gamma(K)}e^{-x(1+\gamma')/\gamma'}\dee x \\
\le & \gamma'^{2(2-\lambda)}\int_{0}^{\gamma'}\frac{x^{K-1+2(\lambda-2)}}{\Gamma(K)}\dee x + (1-e^{-1})^{\lambda-2}\int_{0}^{\infty}\frac{x^{K-1}}{\Gamma(K)}e^{-x(1+\gamma')/\gamma'}\dee x  \\
=& \frac{1}{\Gamma(K)}\frac{\gamma'^K}{K + 2(\lambda-2)} + (1-e^{-1})^{\lambda-2}\left(\frac{\gamma'}{1+\gamma'}\right)^K.
\]
As $K\rightarrow\infty$, the second term will dominate the first term, so when $\lambda-2<0$, the following inequality holds for large $K$,
\[
B_{K,3}\le 4\gamma(1-e^{-1})^{\lambda-2}\left(\frac{\gamma'}{1+\gamma'}\right)^K. 
\]
$B_{K} = B_{K,1} + B_{K,2} + B_{K,3}$ and as $K\rightarrow\infty$, $B_{K,3}$ will dominate $B_{K,1}$ and $B_{K,2}$. So there exists a $K_0\in \nats$ such that for $K > K_0$,
\[
B_{K} \le 12\gamma(1-e^{-1})^{\lambda-2}\left(\frac{\gamma'}{1+\gamma'}\right)^K\longrightarrow 0. 
\]
\paragraph{Sparse network}
When $\alpha>0$,
\[
\nu(\dee\theta) = \gamma'\theta^{-1-\alpha}(1-\theta)^{\lambda+\alpha-1}\dee\theta,\qquad \mu(\dee\theta) = \gamma'\theta^{-1-\alpha} \dee\theta,
\]
and
\[
\frac{\dee\nu}{\dee\mu} = (1-\theta)^{\lambda+\alpha-1}, & & \mu\left([x, 1]\right) = \gamma'\alpha^{-1}(x^{-\alpha} - 1),& & \muinv(u) = \left(1+\frac{\alpha u}{\gamma'}\right)^{-1/\alpha}. 
\]
Similar to the case when $\alpha=0$, 
\[
B_{K,1} = 2\int_{0}^{1}\int_{x_1}^{1}-\log(1-x_1x_2)F_K\left(\mu([x_2, 1])\right)\gamma' x_2^{-1-\alpha}(1-x_2)^{\lambda+\alpha-1}\dee x_2\nu(\dee x_1).
\]
Since $\log(1-x_1x_2)\ge -x_1x_2/(1-x_1)$, 
\[
B_{K,1} &\le 2\gamma'^2\int_{0}^{1}F_K\left(\gamma'\alpha^{-1}(x^{-\alpha}-1)\right)x^{-2\alpha}(1-x)^{2[(\lambda+\alpha)-1]}\dee x \\
&\le 2\gamma'^2\int_{0}^{1}\frac{(\gamma'\alpha^{-1}x^{-\alpha}(1-x^{\alpha}))^K}{\Gamma(K+1)}x^{-2\alpha}(1-x)^{2[(\lambda+\alpha)-1]}\dee x\\
&\le2\gamma'^2\int_{0}^{1}\frac{(\gamma'\alpha^{-1}x^{-\alpha})^K}{\Gamma(K+1)}x^{-2\alpha}(1-x)^{2[(\lambda+\alpha)-1]}\dee x \\
&\le 2\alpha^2\int_{0}^{1}\frac{\left(\gamma'\alpha^{-1}x^{-\alpha}\right)^{K+2}}{{\Gamma(K+1)}}\dee x.
\]
Denoting $t=\gamma'\alpha^{-1}x^{-\alpha}$, then 
\[
B_{K,1} &\le 2\alpha(\gamma'\alpha^{-1})^{1/\alpha} \int_{\gamma'\alpha^{-1}}^{\infty}\frac{t^{K+1-1/\alpha}}{\Gamma(K+1)}\ \dee t \\
&\le 2\alpha(\gamma'\alpha^{-1})^{1/\alpha}\  e^{\gamma'\alpha^{-1}}\ \frac{\Gamma(K+2-1/\alpha)}{\Gamma(K+1)} \int_{\gamma'\alpha^{-1}}^{\infty}\frac{t^{K+1-1/\alpha}}{\Gamma(K+2-1/\alpha)}e^{-t}\ \dee t\\
&\le 2\alpha(\gamma'\alpha^{-1})^{1/\alpha}\  e^{\gamma'\alpha^{-1}}\ \frac{\Gamma(K+2-1/\alpha)}{\Gamma(K+1)}.
\]
By Stirling's formula, 
\[
B_{K,1} \le 2\alpha(\gamma'\alpha^{-1})^{1/\alpha}\  e^{\gamma'\alpha^{-1}}\ K^{\frac{\alpha-1}{\alpha}}.
\]
Now we consider the error bound for $B_{K,2}$. As we have shown in the example when $\alpha=0$, here we can obtain that
\[
B_{K,2} \le 2\alpha(\gamma'\alpha^{-1})^{1/\alpha}\  e^{\gamma'\alpha^{-1}}\ (K-1)^{\frac{\alpha-1}{\alpha}}.
\]
Similar to $B_{K,1}$ and $B_{K,2}$, 
\[
&B_{K,3}\\ 
=& 2\EE\left[\int_{0}^{1}-\log\left(1-\muinv(\Gamma_K)x\right)\frac{\dee \nu}{\dee \mu}\left(\muinv(\Gamma_K)\right)\mathds{1}\left(\Gamma_K\le \muinv[x, 1] \right)\nu( \dee x) \right]\\
\le& 2\gamma'\EE\left[\left(1-\muinv(\Gamma_K)\right)^{\lambda+\alpha-1}\int_{0}^{\muinv(\Gamma_K)}\frac{\muinv(\Gamma_K)x}{1-x}x^{-1-\alpha}(1-x)^{\lambda+\alpha-1}\dee x \right] \\
=& 2\gamma'\EE\left[\left(1-\muinv(\Gamma_K)\right)^{\lambda+\alpha-1}\muinv(\Gamma_K) \int_{0}^{\muinv(\Gamma_K)}x^{-\alpha}(1-x)^{\lambda+\alpha-2}\dee x \right].
\]
We again split our analysis into two cases. First, suppose that $\lambda+\alpha-2\ge0$. Then
\[
&B_{K,3}\\
\le& \frac{2\gamma'}{1-\alpha}\EE\left[\left(1-\muinv(\Gamma_K)\right)^{\lambda+\alpha-1}\muinv(\Gamma_K)^{2-\alpha} \right] \\
=& \frac{2\gamma'}{1-\alpha}\int_{0}^{\infty}\left[1-\left(1+\frac{\alpha x}{\gamma'}\right)^{-1/\alpha}\right]^{\lambda+\alpha-1}\left(1+\frac{\alpha x}{\gamma'}\right)^{-(2-\alpha)/\alpha}\ \frac{x^{K-1}}{\Gamma(K)}e^{-x}\dee x.
\]
Since $\alpha < 1$,
\[
\left(1+\frac{\alpha x}{\gamma'}\right)^{-(2-\alpha)/\alpha}\le \left(\frac{\alpha x}{\gamma'}\right)^{-(2-\alpha)/\alpha}, 
\]
and so
\[
B_{K,3}&\le \frac{2\gamma'}{1-\alpha}\int_{0}^{\infty}\left(\frac{\alpha }{\gamma'}\right)^{-(2-\alpha)/\alpha}\frac{x^{K-1-(2-\alpha)/\alpha}}{\Gamma(K)}e^{-x}\dee x \\
&\le \frac{2\gamma'}{1-\alpha} \left(\frac{\alpha}{\gamma'}\right)^{-(2-\alpha)/\alpha}\ \frac{\Gamma(K-(2-\alpha)/\alpha)}{\Gamma(K)} \\
&\sim \frac{2\gamma'}{1-\alpha} \left(\frac{\alpha}{\gamma'}\right)^{-(2-\alpha)/\alpha} K^{-(2-\alpha)/\alpha},
\]
where the last equation is obtained from Stirling's formula.

On the other hand, if $\lambda + \alpha -2 <0$,
\[
&B_{K,3}\\
\le& 2\gamma'\EE\left[\left(1-\muinv(\Gamma_K)\right)^{\lambda+\alpha-1}\muinv(\Gamma_K) \int_{0}^{\muinv(\Gamma_K)}x^{-\alpha}(1-x)^{\lambda+\alpha-2}\dee x \right]\\
\le& \frac{2\gamma'}{1-\alpha}\EE\left[\left(1-\muinv(\Gamma_K)\right)^{2(\lambda+\alpha)-3}\muinv(\Gamma_K)^{2-\alpha} \right]\\
=& \frac{2\gamma'}{1-\alpha}\int_{0}^{\infty}\left[1-\left(1+\frac{\alpha x}{\gamma'}\right)^{-1/\alpha}\right]^{2(\lambda+\alpha)-3}\left(1+\frac{\alpha x}{\gamma'}\right)^{-(2-\alpha)/\alpha}\ \frac{x^{K-1}}{\Gamma(K)}e^{-x}\dee x.
\]
If $2(\lambda+\alpha)-3\ge 0$, we get the same result as in the case when $\lambda+\alpha-2\ge0$. When $2(\lambda+\alpha)-3\le 0$, note that we can find an $x_0$ such that when $x\in[0, x_0]$, 
\[1-\left(1+\frac{\alpha x}{\gamma'}\right)^{-1/\alpha}\ge x^2. \]
So
\[
&B_{K,3}\\
\le& \frac{2\gamma'}{1-\alpha}\left\{ \int_{0}^{x_0} x^{4(\lambda+\alpha)-6}\left(1+\frac{\alpha x}{\gamma'}\right)^{-(2-\alpha)/\alpha}\ \frac{x^{K-1}}{\Gamma(K)}e^{-x}\dee x\right.\\
&\left. +\left[1-\left(1+\frac{\alpha x_0}{\gamma'}\right)^{-1/\alpha}\right]^{2(\lambda+\alpha)-3}\int_{x_0}^{\infty}\left(1+\frac{\alpha x}{\gamma'}\right)^{-(2-\alpha)/\alpha}\ \frac{x^{K-1}}{\Gamma(K)}e^{-x}\dee x \right\} \\
\le& \frac{\gamma'}{1-\alpha}\left(\frac{\alpha}{\gamma'}\right)^{-(2-\alpha)/\alpha}\left\{\int_{0}^{x_0}\frac{x^{K-1+4(\lambda+\alpha)-6-(2-\alpha)/\alpha}}{\Gamma(K)}e^{-x}\dee x \right.\\
& \left.+ \left[1-\left(1+\frac{\alpha x_0}{\gamma'}\right)^{-1/\alpha}\right]^{2(\lambda+\alpha)-3}\frac{\Gamma(K-(2-\alpha)/\alpha)}{\Gamma(K)} \right\}.
\]
Because here we assume that $2(\lambda + \alpha)-3<0$, the second term will dominate. So in this case we obtain that for large $K$,
\[
B_{K,3}&\le \frac{4\gamma'}{1-\alpha} \left(\frac{\alpha}{\gamma'}\right)^{-(2-\alpha)/\alpha}\left[1-\left(1+\frac{\alpha x_0}{\gamma'}\right)^{-1/\alpha}\right]^{2(\lambda+\alpha)-3} K^{-(2-\alpha)/\alpha}.
\]
Asymptotically, $B_{K,2}$ will dominate $B_{K,1}$ and $B_{K,3}$, so there exists $K_0\in \nats$ such that when $K > K_0$,
\[
B_{K}\le 6\alpha(\gamma'\alpha^{-1})^{1/\alpha}\  e^{\gamma'\alpha^{-1}}\ (K-1)^{\frac{\alpha-1}{\alpha}} \longrightarrow 0.
\]
\subsection{Gamma-independent Poisson network}\label{pf:gamma}
\paragraph{Dense network} When $\alpha = 0$,
\[
\nu(\dee\theta) = \gamma\lambda\theta^{-1}e^{-\lambda\theta}\dee \theta, \qquad \mu(\dee \theta) = \gamma\lambda\theta^{-1}(1+\lambda\theta)^{-1}\dee\theta. 
\]
In this case,
\[
\frac{\dee\nu}{\dee\mu} = (1+\lambda\theta)e^{-\lambda\theta}, & & \mu[x, \infty) = \gamma\lambda\log(1+(\lambda x)^{-1}),& & \muinv(u) = \frac{1}{\lambda(e^{(\gamma\lambda)^{-1}u} - 1)}. 
\]
For Poisson distribution, $\pi(\theta) = e^{-\theta}$, so 
\[
B_{K,1} = \int_{\reals_+^2}x_1x_2F_K\left(\mu[\max\{x_1, x_2 \}, \infty) \right)\nu(\dee x_1)\nu(\dee x_2).
\]
Note that the integrand is symmetric around the line $x_1 = x_2$, so we only need to compute the integral above the line $x_1 = x_2$. In this region, $e^{-\lambda x_2} \le e^{-\lambda x_1}$ and $0\le F_K\left(\mu[x_2, \infty)\right)\le F_K\left(\mu[x_1, \infty)\right)\le 1$. So
\[
B_{K,1} &= 2\int_{\reals_+}x_1\int_{x_1}^{\infty}x_2F_K\left(\mu[x_2, \infty) \right)\lambda\gamma x_2^{-1}e^{-\lambda x_2} \dee x_2 \nu(\dee x_1) \\
&\le 2\int_{\reals_+}x_1 F_K\left(\mu[x_1, \infty) \right)\gamma e^{-\lambda x_1}\nu(\dee x_1) \\
&= 2\gamma^2\lambda\int_{\reals_+}F_K\left(\mu[x, \infty) \right)e^{-2\lambda x} \dee x.
\]
For any $a>0$, we divide the integral into two parts and bound each part separately. We denote $b = \log(1+(\lambda a)^{-1})$ and use the fact that $\int_0^a e^{-2\lambda x}\dee x \le a$ and  $F_K(t)\le t^K/K!\le(3t/K)^K$. So
\[
B_{K,1}&\le 2\gamma^2\lambda \left[\int_{0}^{a}e^{-2\lambda x}\dee x + F_K\left(\mu[a, \infty) \right)\int_{a}^{\infty}e^{-2\lambda x}\dee x \right]\\
&\le 2\gamma^2\lambda\left[a + \frac{1}{2\lambda}F_K\left(\gamma\lambda\log(1+(\lambda a)^{-1}) \right) \right] \\
&\le 2\gamma^2\lambda\left[\left(\lambda(e^b-1) \right)^{-1} + \frac{1}{\lambda}\left(\frac{3\gamma\lambda b}{K}\right)^K \right].
\]
Setting two terms equal and use the fact that $(e^b -1)^{-1} \approx e^{-b}$ when $b$ is large, we get $b = KW_0((3\gamma\lambda)^{-1})$ and $W_0$ is defined by
\[
W_0(y) = x \iff xe^x=y. 
\]
Therefore,
\[
B_{K,1}\le \frac{4\gamma^2}{e^{KW_0((3\gamma\lambda)^{-1})}-1} \sim 4\gamma^2e^{-KW_0\left((3\gamma\lambda)^{-1}\right)}. 
\]
Similarly, 
\[
&\frac{1}{2(K-1)}B_{K,2} \\
=&\int_{\reals_+^2}x_1x_2\EE\left[\frac{1}{\Gamma_K}\mathds{1}\left(\mu[x_2, \infty)\le \Gamma_K \le \mu[x_1, \infty) \right) \right]\nu(\dee x_1)\nu(\dee x_2). 
\]
Note that
\[
&\EE\left[\frac{1}{\Gamma_K}\mathds{1}\left(\mu[x_2, \infty)\le \Gamma_K \le \mu[x_1, \infty) \right) \right] \\
=& \frac{1}{K-1}\left[F_{K-1}\left(\mu[x_1, \infty)\right)-F_{K-1}\left(\mu[x_2, \infty) \right) \right].
\]
Keeping only the positive part,
\[
B_{K,2} &\le 2\int_{\reals_+}x_1 F_{K-1}\left(\mu[x_1, \infty) \right)\int_{x_1}^{\infty}\gamma\lambda e^{-\lambda x_2}\dee x_2\nu(\dee x_1) \\
&=2\gamma^2\lambda\int_{\reals_+}F_{K-1}\left(\mu[x,\infty) \right)e^{-2\lambda x} \dee x.
\]
This has the same form as $B_{K,1}$, so 
\[
B_{K,2}\le \frac{4\gamma^2}{e^{(K-1)W_0((3\gamma\lambda)^{-1})}-1} \sim 4\gamma^2e^{-(K-1)W_0\left((3\gamma\lambda)^{-1}\right)}. 
\]
Next,
\[
&B_{K,3} \\
=& 2\EE\left[\int_{\reals_+}x\muinv(\Gamma_K)\frac{\dee \nu}{\dee\mu}\left(\muinv(\Gamma_K) \right)\mathds{1}\left[\Gamma_K\le\mu [x, \infty) \right] \right]\nu(\dee x) \\
=& 2\EE\left[\int_{0}^{\muinv(\Gamma_K)}x\muinv(\Gamma_K)\left(1+\lambda\muinv(\Gamma_K)\right)e^{-\lambda\muinv(\Gamma_K)}\gamma\lambda x^{-1}e^{-\lambda x}\dee x \right] \\
=&2\gamma\EE\left[\muinv(\Gamma_K)\left(1+\lambda\muinv(\Gamma_K)\right)e^{-\lambda\muinv(\Gamma_K)}\left(1-e^{-\lambda\muinv(\Gamma_K)} \right) \right] \\
\le& 2\gamma\lambda\EE\left[\muinv(\Gamma_K)^2\left(1+\lambda\muinv(\Gamma_K)\right)e^{-\lambda\muinv(\Gamma_K)} \right].
\]
Note that $(1+x)e^{-x} \le 1$, so
\[
B_{K,3} &\le 2\lambda\gamma \EE\left[\muinv(\Gamma_K)^2\right] \\
&\le 2\frac{\gamma}{\lambda}\EE\left[e^{-(\gamma\lambda)^{-1}\Gamma_K} \right] 
= 2\frac{\gamma}{\lambda}\left(\frac{\gamma\lambda}{1+\gamma\lambda}\right)^{K-1}.
\]
Since $B_{K,3}$ will dominate $B_{K,1}$ and $B_{K,2}$ asymptotically, there exists $K_0 \in \nats$ such that for $K > K_0$,
\[
B_{K} \le 3B_{K,3}\le 6\frac{\gamma}{\lambda}\left(\frac{\gamma\lambda}{1+\gamma\lambda}\right)^{K-1}. \]

\paragraph{Sparse network}
When $\alpha > 0$,
\[
\nu(\dee\theta) = \gamma\frac{\lambda^{1-\alpha}}{\Gamma(1-\alpha)}\theta^{-\alpha-1}e^{-\lambda\theta}\dee\theta, \qquad \mu(\dee\theta) = \gamma\frac{\lambda^{1-\alpha}}{\Gamma(1-\alpha)}\theta^{-\alpha-1}\dee\theta, 
\]
and
\[
\frac{\dee\nu}{\dee\mu} = e^{-\lambda\theta},& & \mu[x, \infty) = \gamma'x^{-\alpha}, & & \muinv(u)=(\gamma'u^{-1})^{1/\alpha},& & \gamma' = \gamma\frac{\lambda^{1-\alpha}}{\alpha\Gamma(1-\alpha)}. 
\]
Similar to the example when $\alpha=0$,
\[
B_{K,1} &= \int_{\reals_+^2}x_1x_2F_K\left(\mu[\max\{x_1, x_2\}, \infty) \right)\nu(\dee x_1)\nu(\dee x_2) \\
&=  2\int_{\reals_+}x_1\int_{x_1}^{\infty}F_K\left(\mu[x_2, \infty) \right)\gamma\frac{\lambda^{1-\alpha}}{\Gamma(1-\alpha)}x_2^{-\alpha}e^{-\lambda x_2}\dee x_2 \nu(\dee x_1) \\
&\le 2\gamma\int_{\reals_+}x_1F_K\left(\mu[x_1, \infty) \right)\int_{x_1}^{\infty}\frac{\lambda^{1-\alpha}}{\Gamma(1-\alpha)}x_2^{-\alpha}e^{-\lambda x_2}\dee x_2 \nu(\dee x_1).
\]
Note that the integrand with respect to $x_2$ is the density function of the
gamma distribution with shape $\alpha$ and rate $\lambda$, so the integral is less than 1.
We partition the outer integral into two parts and bound them separately, 
\[
B_{K,1}&\le 2\gamma\int_{\reals_+}F_K\left(\mu[x, \infty)\right)\gamma\frac{\lambda^{1-\alpha}}{\Gamma(1-\alpha)}x^{-\alpha}e^{-\lambda x}\dee x  \\
&\le 2\gamma^2\frac{\lambda^{1-\alpha}}{\Gamma(1-\alpha)}\left[\int_{0}^{a}x^{-\alpha}\dee x + F_K\left(\mu[a, \infty) \right)\int_{a}^{\infty}x^{-\alpha}e^{\lambda x}\dee x \right] \\
&\le 2\gamma^2\frac{1}{\Gamma(1-\alpha)}\left[\frac{\lambda^{1-\alpha}}{1-\alpha}a^{1-\alpha} + \Gamma(1-\alpha)\left(\frac{3\gamma'a^{-\alpha}}{K} \right)^K \right].
\]
By setting the two terms in the brackets equal, we get 
\[
a = \left[\frac{(1-\alpha)\Gamma(1-\alpha)}{\lambda^{1-\alpha}} \right]^{\frac{1}{(K-1)\alpha +1}}\left(\frac{3\gamma'}{K}\right)^{\frac{K}{(K-1)\alpha +1}}\sim \left(\frac{3\gamma'}{K}\right)^{\frac{1}{\alpha}}. 
\]
So
\[
B_{K,1} \le \frac{4\gamma^2\lambda^{1-\alpha}}{(1-\alpha)\Gamma(1-\alpha)}\left(\frac{3\gamma'}{K}\right)^{\frac{1-\alpha}{\alpha}}. 
\]
Similar to the last example where $\alpha=0$, here 
\[
&B_{K,2} \\
=& 2(K-1)\int_{\reals_+^2}x_1x_2\EE\left[\frac{1}{\Gamma_K}\mathds{1}\left(\mu[x_2, \infty)\le \Gamma_K \le \mu[x_1,\infty) \right) \right]\nu(\dee x_1)\nu(\dee x_2) \\
\le& 2\int_{\reals_+}x_1F_{K-1}\left(\mu[x_1, \infty) \right)\int_{x_1}^{\infty}x_2\nu(\dee x_2)\nu(\dee x_1) \\
\le& 2\gamma\int_{\reals_+}x_1F_{K-1}\left(\mu[x_1, \infty) \right)\nu(\dee x_1).
\]
This has the same form as $B_{K,1}$, and therefore
\[
B_{K,2}\le \frac{4\gamma^2\lambda^{1-\alpha}}{(1-\alpha)\Gamma(1-\alpha)}\left(\frac{3\gamma'}{K-1}\right)^{\frac{1-\alpha}{\alpha}}. 
\]
Finally, since both $e^{-\lambda x} <1$ and $e^{-\lambda\muinv(\Gamma_K)} < 1$, 
\[
&B_{K,3}\\
=& 2\EE\left[\int_{\reals_+}\muinv(\Gamma_K)x e^{-\lambda\muinv(\Gamma_K)}\mathds{1}\left[x\le\muinv(\Gamma_K) \right]\gamma\frac{\lambda^{1-\alpha}}{\Gamma(1-\alpha)}x^{-1-\alpha}e^{-\lambda x}\dee x \right] \\
\le& 2\gamma\frac{\lambda^{1-\alpha}}{\Gamma(1-\alpha)}\EE\left[\muinv(\Gamma_K) \int_{0}^{\muinv(\Gamma_K)}x^{-\alpha}\dee x \right] \\
\le& 2\gamma\frac{\lambda^{1-\alpha}}{(1-\alpha)\Gamma(1-\alpha)}\EE\left[\muinv(\Gamma_K)^{2-\alpha} \right] \\
=& 2\gamma\frac{\lambda^{1-\alpha}}{(1-\alpha)\Gamma(1-\alpha)}\left(\gamma'\right)^{\frac{2-\alpha}{\alpha}}\frac{\Gamma(K-\frac{2-\alpha}{\alpha})}{\Gamma(K)}.
\]
By Stirling's formula,
\[
\frac{\Gamma(K-(2-\alpha)/\alpha)}{\Gamma(K)}\sim \frac{\sqrt{2\pi(K-\frac{2-\alpha}{\alpha})}\left(\frac{K-(2-\alpha)/\alpha}{e}\right)^{K-(2-\alpha)/\alpha} }{\sqrt{2\pi K}\left(\frac{K}{e}\right)^K }\sim K^{-\frac{2-\alpha}{\alpha}}. 
\]
So
\[
B_{K,3} \le 2\gamma\frac{\lambda^{1-\alpha}}{(1-\alpha)\Gamma(1-\alpha)}\left(\gamma'\right)^{\frac{2-\alpha}{\alpha}} K^{-\frac{2-\alpha}{\alpha}}. 
\]
$B_{K,2}$ dominates $B_{K,1}$ and $B_{K,3}$ asymptotically, so there exists $K_0\in \nats$ such that for $K > K_0$,
\[
B_{K}\le \frac{12\gamma^2\lambda^{1-\alpha}}{(1-\alpha)\Gamma(1-\alpha)}\left(\frac{3\gamma'}{K-1}\right)^{\frac{1-\alpha}{\alpha}}. 
\]

\section{Truncated inference}\label{sec:trunc_gibbs_derivation}

\bprfof{\cref{thm:mherror}}
Fix $K\in\nats$ and $\epsilon >0$, and define the 
subset of state space $A = \left\{ X_{N+1}^{K+} = 0, B(\Gamma_{1:K}, \sigma) \le \frac{\epsilon}{N+1}\right\}$. 
By assumption, 
\[
\hPi\left(A\right) \geq 1-\eta.
\]
Further, by applying the bound from \cref{thm:conditionaltrunc}, we know that for states in $A$,
\[
1-\epsilon \leq p(X_{1:N+1}^{K+} \given \theta_{1:K}, \sigma) \leq 1.
\]
Suppose $p$ is the RHS of \cref{eq:postdensity} that is proportional to the density of $\Pi$, and $\hp$ is the RHS of \cref{eq:postdensity} removing the term $p(X_{1:N}^{K+} | \theta_{1:K}, \sigma)$, which is proportional to the density of $\hat{\Pi}$:
\[
p \leq \hp, \qquad \text{and within $A$,}\quad p \geq (1-\epsilon)\hp.
\]
Define the normalization constants $Z, \hZ$ such that $\int p/Z = \int \hp/\hZ = 1$.
Then the above bounds yield
\[
\hZ = \!\int \hp \geq \!\int p
 = \!\int_A p + \!\int_{A^c} p
 \geq \!\int_A(1-\epsilon)\hp
 = (1-\epsilon)\frac{\int_A\hp}{\int \hp} \! \int \hp
 \geq (1-\epsilon)(1-\eta) \hZ.
\]
Therefore, $(1-\epsilon)(1-\eta)\hZ \leq Z \leq \hZ$, and hence
\[
\Pi(A) = \int_A p/Z \geq (1-\epsilon) \int_A \hp/Z
\geq (1-\epsilon) \int_A \hp/\hZ
\geq (1-\epsilon)(1-\eta).
\]
The above results yield the total variation bound via
\[
\frac{1}{2} \int \left| \frac{p}{Z} - \frac{\hp}{\hZ}\right| &=
\frac{1}{2} \int_A \left| \frac{p}{Z} - \frac{\hp}{\hZ}\right|  +
\frac{1}{2} \int_{A^c} \left| \frac{p}{Z} - \frac{\hp}{\hZ}\right| \\
&\leq
\frac{1}{2} \int_A \left| \frac{p}{Z} - \frac{p}{\hZ}\right|  +
\frac{1}{2}\int_A \left| \frac{p}{\hZ} - \frac{\hp}{\hZ}\right|  +
\frac{1}{2} \left( \Pi(A^c) + \hPi(A^c)\right)\\
&=
\frac{1}{2} \left| \frac{1}{Z} - \frac{1}{\hZ}\right|\int_A p  +
\frac{1}{2}\frac{1}{\hZ} \int_A (\hp-p)  +
\frac{1}{2} \left( \Pi(A^c) + \hPi(A^c)\right)\\
&\leq
\frac{1}{2} \left( 1 - (1-\epsilon)(1-\eta)\right)  +
\frac{1}{2}\epsilon   +
\frac{1}{2} \left(1 - (1-\eta)(1-\epsilon) + \eta\right)\\
&=\frac{3(\epsilon+\eta)}{2} - \epsilon\eta.
\]
\eprfof
\end{document}